\documentclass{article}
\usepackage[english]{babel}
\usepackage[utf8]{inputenc}
\usepackage{array} 
\usepackage{multirow} 
\usepackage{amsmath, amssymb, graphicx}
\usepackage{caption}
\usepackage{bm}
\usepackage{graphicx}
\usepackage{color}
\usepackage{multirow}
\usepackage{algorithm}
\usepackage{algpseudocode}
\usepackage{subcaption}
\usepackage{amsfonts}
\usepackage{tikz}
\usepackage{placeins}
\usepackage{todonotes}
\usetikzlibrary{arrows.meta,positioning,fit,calc}

\usepackage{comment}

\graphicspath{ {./figs/} }
\newcommand{\hr}[1]{{\color{black}#1}}
\newcommand{\hb}[1]{{\color{black}#1}}
\newcommand{\rev}[1]{{\color{black}#1}}

\title{
\rev{Learning Hidden Ecological and Epidemic Processes\\ in a Reduced-Order Space}
}

\author{
Maria Vasilyeva 
\thanks{Department of Mathematics and Statistics, Texas A\&M University-Corpus Christi,   Corpus Christi, TX, USA\\
Email: {\tt maria.vasilyeva@tamucc.edu}.} 
\and
Zheng Wei\footnotemark[1] 
\and 
Kelum Gajamannage\thanks{Department of Mathematics and Applied Mathematical Sciences, University of Rhode Island, Kingston, RI, USA\\
Email: {\tt kelum.gajamannage@uri.edu}.} 
\and
Hyangim Ji\footnotemark[1] 
\and
Aleksei Krasnikov\footnotemark[1] 
\and 
Alexey Sadovski\footnotemark[1] 
}

\begin{document}

\maketitle

\begin{abstract}We study data-driven learning of ecological and epidemic dynamics governed by nonlinear reaction and reaction--diffusion systems. We consider both neural approaches for learning the underlying reaction dynamics from trajectory data and Sparse Identification of Nonlinear Dynamics (SINDy) for discovering sparse, interpretable governing equations.
A key component of this work is a partial-learning framework for reaction--diffusion systems. The prescribed diffusion operator is projected onto a reduced-order space and explicitly retained during differentiable integration, while only the unknown local eco-epidemic reaction is learned. This enables efficient training on fine spatial grids by evaluating the neural reaction only at a small set of informative nodes. The reduced solution is then reconstructed on the original fine grid and evaluated over both the training interval and an unseen future prediction horizon.
Finally, the learned neural reaction is processed by SINDy to recover a sparse, interpretable symbolic model of the underlying dynamics. Experiments with clean and noisy data demonstrate high predictive accuracy, implicit denoising, computational efficiency, and the ability to recover sparse reaction models.
\end{abstract}

\section{Introduction}

Ecological and epidemic systems play a crucial role in shaping population survival, maintaining biodiversity, and ensuring ecosystem stability  \cite{prenter2004roles, anderson1991populations, su2015modeling}. \rev{When ecological and epidemic processes are coupled, they produce nonlinear eco--epidemic dynamics governed by trophic interactions and disease transmission, while spatial dispersal introduces interactions across locations~\cite{ji2025epidemic}. In this work, we present a hybrid data-driven framework for reduced-order learning of spatial eco--epidemic systems. The diffusion operator is retained explicitly in the reduced space, while only the unknown local reaction dynamics are learned from data. }

Recent advances in data-driven modeling have enabled the discovery of governing equations directly from time-series data \cite{weinan2017proposal, park2022recurrent}. Several neural network architectures have been developed to complement or replace classical equation-based approaches, particularly in the context of epidemic modeling. Neural networks (NNs) provide flexible frameworks capable of capturing complex nonlinear dynamics commonly observed in epidemics \cite{datilo2019review}. Neural ODEs have been introduced for modeling epidemic dynamics, providing a continuous-time framework that naturally aligns with the underlying differential equation structure of epidemic processes  \cite{kosma2023neural, huang2025neuralcode}. Epi-DNN employs deep neural networks to estimate time-varying or unknown parameters of traditional epidemic models, such as transmission and recovery rates \cite{ning2023epi}. The Long Short-Term Memory (LSTM) algorithm, which uses a type of recurrent neural network (RNN), was previously used in other studies to predict infections over time \cite{zhu2019attention, chae2018predicting}. \rev{Physics-Informed Neural Networks (PINNs) provide another important framework by enforcing a prescribed epidemiological residual during forward or inverse learning~\cite{raissi2019physics,kharazmi2021identifiability,millevoi2024physics}.}
\rev{Neural networks} have been successfully applied in ecology to model species distributions and predict ecological outcomes. Artificial neural networks improve predictive performance across diverse ecological contexts \cite{lek1999artificial}. 
In predator–prey systems, continuous-time neural networks can simulate dynamic ecological interactions \cite{moreau1999embedding}, and training algorithms tailored to predator–prey dynamics help optimize network parameters for modeling complex dynamical behaviors \cite{tilahun2013prey}. \rev{Nevertheless, ecological learning and epidemic forecasting have largely developed as separate strands, and fewer studies address interpretable learning of their coupled dynamics.}

\rev{Resolving spatial interactions introduces a second challenge: the learning method must represent dependence across locations while evolving the system in time. For spatially resolved data, convolutional neural networks exploit the regular neighborhood structure of images and structured grids~\cite{lecun1998gradient,ronneberger2015unet,xiao2024numerical}, whereas graph neural networks extend neighborhood-based learning to unstructured meshes and interaction graphs~\cite{scarselli2008graph,kipf2017semi}. Spatio-Temporal Graph Neural Networks (STGNNs) represent locations or subpopulations as graph nodes and combine graph-based spatial aggregation with recurrent, convolutional, or attention-based temporal encoders; applications include regional epidemic transmission and metapopulation forecasting~\cite{xie2022epignn,cao2022mepognn,jin2024stgnnsurvey}. Continuous-time variants replace discrete latent updates with graph neural ODEs that evolve spatial and temporal representations continuously~\cite{fang2021stgode,jin2023mtgode}. 
These methods are natural spatial forecasting baselines when mobility or contact networks define the interactions. 
In this work, space is represented by a fine-grid discretization of a reaction--diffusion system, where the prescribed diffusion operator is evaluated explicitly throughout the learning process, while the local eco--epidemic reaction is the unknown quantity to be learned.}

\rev{Fine-grid spatial discretizations can produce large systems, making repeated time integration and gradient-based training computationally expensive.} \rev{Proper Orthogonal Decomposition (POD) reduces this cost by extracting a low-dimensional basis from solution snapshots and projecting the state and spatial operators onto that basis~\cite{sirovich1987turbulence,chatterjee2000introduction,vasilyeva2026operator}. Projection alone, however, does not remove the need to reconstruct and evaluate a nonlinear reaction over the full grid. The Discrete Empirical Interpolation Method (DEIM) addresses this bottleneck by approximating the nonlinear term from evaluations at selected spatial nodes~\cite{chaturantabut2010nonlinear}, while QDEIM provides a stable QR-based strategy for choosing those nodes~\cite{drmac2016new}. Related biological applications support the relevance of reduced coordinates: POD with sparse nonlinear sampling has been applied to spatial predator--prey communities~\cite{dimitriu2017comparative}; POD bases built from parametric SIR solutions have supported epidemic data assimilation and forecasting~\cite{chacon2020predictive}; and POD-based compression has been used to learn the evolution of spatial SEIRS models~\cite{silva2023predgan}. These studies demonstrate that ecological and epidemic trajectories can admit informative low-dimensional representations, but they primarily reduce the complete model or learn a surrogate for the reduced state evolution.}

\rev{In this work, we investigate data-driven learning of epidemic and ecological dynamics in temporal and spatial settings. For temporal systems, the complete vector field is inferred from trajectory data, whereas for spatial reaction--diffusion systems, the prescribed diffusion operator is evaluated explicitly and only the unknown local reaction is learned.} Using synthetic datasets \rev{from} canonical ecological and epidemiological models, we \rev{compare} three data-driven frameworks: Sparse Identification of Nonlinear Dynamics (SINDy), Neural Ordinary Differential Equations (Neural ODEs), and Kolmogorov--Arnold Network ODEs (KANODEs). \rev{SINDy seeks parsimonious governing equations through sparse regression~\cite{brunton2016discovering}, whereas Neural ODEs learn flexible, nonparametric dynamics directly from observed trajectories~\cite{chen2018neural}. KANODEs use the Kolmogorov--Arnold network architecture~\cite{liu2024kan,koenig2024kan} to represent nonlinear dependencies with relatively few parameters. Together, these non-spatial comparisons establish the data-driven tools used next for the spatial learning problem. For spatial reaction--diffusion systems, we consider the partial-learning setting in which the diffusion operator is known and only the hidden local reaction is learned. Our central contribution is to formulate this problem in a reduced-order space: POD reduces the evolving state and the prescribed diffusion operator, while DEIM hyper-reduces only the learned nonlinear reaction. The methodological novelty lies specifically in adapting POD--QDEIM to the known/unknown decomposition of physics-structured partial learning: rather than reducing the complete vector field, we retain the prescribed diffusion as a fixed projected operator and construct the nonlinear basis and interpolation nodes from diffusion-subtracted reaction snapshots, so that only the unknown neural reaction is hyper-reduced. The reaction-informed construction permits neural-reaction evaluation at only a small set of spatial nodes, preserves the known diffusion explicitly, and reconstructs predictions on the original fine grid. The reduced-space formulation also decreases training time while maintaining comparable full-trajectory accuracy.}

The paper is organized as follows. \rev{Section~2 presents the temporal benchmark models and the corresponding NODE, KANODE, and SINDy learning formulations. Section~3 develops the spatial partial-learning framework for hidden eco--epidemic coupling in full and reduced spaces. Section~4 presents the numerical experiments, and the paper ends with conclusions in Section~5.}

\section{\rev{Temporal Models and Data-Driven Learning}}

\rev{We first consider temporal ecological and epidemic dynamics governed by ODE models. The SIR system, an LV model augmented with density-dependent prey mortality, and the LVSIS system represent epidemic, ecological, and coupled eco--epidemiological behavior, respectively, and are formulated within a common dynamical framework. Their numerical trajectories provide observations for evaluating data-driven recovery of the underlying vector fields. We then discuss Neural ODEs using MLP and KAN parameterizations together with SINDy as data-driven approaches for learning the temporal dynamics.}

\subsection{\rev{Temporal Model Formulation}}

\rev{We consider ecological and epidemic dynamics described by a general system of ordinary differential equations. As representative benchmarks, we use the Susceptible--Infected--Recovered (SIR), a Lotka--Volterra (LV) model augmented with density-dependent prey mortality, and the coupled Lotka--Volterra--SIS (LVSIS) model, which provide a range of nonlinear dynamics for evaluating the proposed learning framework. Spatial diffusion is incorporated later in Section~3, where these models are extended to reaction--diffusion systems.}

We consider the general ODE system 
\begin{equation}
\label{eq:mm}
\frac{d y}{d t} = f(y), 
\qquad y(0) = y_0,
\end{equation}
where $y(t)$ is the state vector and $f(t,y)$ defines the model dynamics.  
Below are several instances of $f(y)$ used in this work \cite{tolles2020modeling, chattopadhyay1999predator,han2001four,singh2018mathematical}.

\begin{itemize}
\item \textit{SIR model}: 
State vector $y=(u_s,u_i,u_r)^\top$.  
\begin{equation}
f_{\mathrm{SIR}}(y) =
\begin{pmatrix}
-\sigma u_s u_i \\[6pt]
\sigma u_s u_i - \omega u_i \\
\omega u_i
\end{pmatrix}, 
\label{sir}
\end{equation}
with $y(0) = (u_{s0}, u_{i0}, u_{r0})^\top$. 

\item \textit{Lotka--Volterra predator--prey model}: 
State vector $y=(u,v)^\top$ with $u$ = prey and $v$ = predator.  
\begin{equation}
f_{\mathrm{LV}}(y) =
\begin{pmatrix}
\alpha u - \beta u v  - \eta u^2 \\
\delta u v - \gamma v
\end{pmatrix}, 
\label{lv}
\end{equation}
with $y(0) = (u_0, v_0)^\top$. 
\rev{The term $-\eta u^2$ represents density--dependent prey mortality.}

\item \textit{Lotka--Volterra--SIS (eco-epidemiological)}: 
State vector $y=(u_s, u_i, v_s, v_i)^\top$ with susceptible/infected prey and predator populations.  
\begin{equation}
f_{\mathrm{LVSIS}}(y) =
\begin{pmatrix}
\alpha_s u_s + \alpha_i u_i
- \beta^s_s u_s v_s - \beta^s_i u_s v_i
- \gamma^u_s u_s^2
- \sigma^s_i u_s u_i
- \sigma^s_{vi} u_s v_i
+ \omega_i u_i
\\
- \beta^i_s u_i v_s - \beta^i_i u_i v_i
- \gamma^u_i u_i^2
+ \sigma^s_i u_s u_i
+ \sigma^s_{vi} u_s v_i
- \omega_i u_i
\\
\quad- \gamma^v_s v_s + \delta^s_s v_s u_s + \delta^s_i v_s u_i
- \sigma^{vs}_i v_s u_i - \sigma^{vs}_{vi} v_s v_i
+ \omega_{vi} v_i
\\
\quad- \gamma^v_i v_i + \delta^i_s v_i u_s + \delta^i_i v_i u_i
+ \sigma^{vs}_i v_s u_i + \sigma^{vs}_{vi} v_s v_i
- \omega_{vi} v_i
\end{pmatrix},
\label{LVSIS}
\end{equation}
with $y(0) = (u_{s0}, u_{i0}, v_{s0}, v_{i0})^\top$.
\end{itemize}

\rev{For the SIR model, $\sigma$ and $\omega$ denote the disease-transmission and recovery rates, respectively. In the LV model, $\alpha$ is the prey growth rate, $\beta$ the predation rate, $\eta$ the density-dependent prey mortality rate, $\delta$ the predator reproduction rate, and $\gamma$ the predator mortality rate. In the LVSIS model, $u$ and $v$ denote prey and predator populations, while $s$ and $i$ denote susceptible and infected classes. The parameter families $\alpha$, $\beta$, $\gamma^u$, $\gamma^v$, $\delta$, $\sigma$, and $\omega$ represent prey growth, predation, prey density-dependent mortality, predator mortality, predator reproduction, disease transmission, and recovery to the susceptible class, respectively. For the interaction coefficients $\beta$, $\delta$, and $\sigma$, the superscript identifies the affected compartment and the subscript identifies the interacting compartment. Thus, $\beta_i^s u_sv_i$ represents predation of susceptible prey $u_s$ by infected predators $v_i$, $\delta_i^s v_su_i$ represents the contribution of infected prey $u_i$ to susceptible-predator $v_s$ reproduction, and $\sigma_{vi}^{vs}v_sv_i$ represents disease transmission from infected predators $v_i$ to susceptible predators $v_s$. The parameter values used in the experiments are specified in Section~4.}

To numerically solve the time-dependent problem \eqref{eq:mm} on the interval $0 \le t \le T$, time integration is performed using an explicit Runge--Kutta (RK) method. Given a time step $\tau_n$ and the numerical solution $y^n \approx y(t^n)$ at time $t^n$, an $s$-stage RK method is defined by
\begin{align}
k_i &= f\!\left(y^n + \tau_n \sum_{j=1}^{i-1} a_{ij} k_j\right), 
\qquad i = 1,\ldots,s, \label{eq:rk-stages} \\
y^{n+1} &= y^n + \tau_n \sum_{i=1}^{s} b_i k_i. \label{eq:rk-update}
\end{align}
The coefficients $\{a_{ij}, b_i\}$ characterize the specific RK scheme and are typically represented in a Butcher tableau. 
In this work, we employ an explicit Runge--Kutta method of order $4$.

\subsection{\rev{Learning Temporal Dynamics}}

Classical modeling frameworks require explicit knowledge of $f(y)$ in \eqref{eq:mm}. However, in many real-world eco-epidemiological systems, the exact functional form of the dynamics is unknown or too complex to be expressed analytically. 
We consider several data-driven approaches, including Neural ODEs (NODEs), Kolmogorov--Arnold Network ODEs (KANODEs), and Sparse Identification of Nonlinear Dynamics (SINDy), to learn the underlying dynamics directly from observational data.  
\rev{Because these methods are well established, we summarize only the formulations needed in this work.}

Given a set of observed data $\{y^n_i\}_{i=1}^N$ at time $t_n$, we define a loss function to quantify the difference between the predicted values and the observations, 
\begin{equation}
\mathcal{L}(\theta) = \frac{1}{N} \sum_{n=1}^N \|y_\theta^n - y^n\|^2,
\label{eq:loss}
\end{equation}
\rev{where $y_\theta^n$ denotes the predicted value at time $t_n$.  
The observations are generated from numerical simulations of the temporal models formulated above.} A fixed time-stepping scheme is employed to generate the time series, to which random noise is added. We then investigate the ability of the considered data-driven approaches to learn the underlying dynamics from these noisy observations.
\rev{Figure~\ref{fig:learning-workflow-ode} summarizes the shared workflow from data generation and preparation to model training, prediction, and evaluation.}

\begin{figure}[h!]
\centering
\resizebox{0.6\textwidth}{!}{%
\begin{tikzpicture}[
    font=\small, >=Latex, node distance=8mm and 9mm,
    flowbox/.style={draw=black,text=black,rounded corners,align=center,
        minimum height=12mm,minimum width=33mm,fill=black!4},
    method/.style={draw=black,text=black,rounded corners,align=center,
        minimum height=18mm,minimum width=38mm,fill=black!8},
    flow/.style={->,thick,draw=black}
]
\node[flowbox] (systems) {SIR, LV, LVSIS\\synthetic trajectories};
\node[flowbox,right=of systems] (prep) {Noise and normalization,\\training interval /\\unseen future horizon};
\draw[flow] (systems)--(prep);

\node[method,below=16mm of prep,xshift=-24mm,minimum width=48mm] (node)
    {\textbf{Neural ODE (NODE)}\\$\dot y=\mathrm{NN}(y;\theta)$\\NN: MLP or KAN};
\node[method,right=18mm of node] (sindy)
    {\textbf{SINDy}\\library $\Theta(Y)$\\sparse regression};

\coordinate (branch) at ($(prep.south)+(0,-6mm)$);
\draw[draw=blue,thick] (prep.south)--(branch);
\draw[flow] (branch)-|(node.north);
\draw[flow] (branch)-|(sindy.north);

\node[flowbox,below=16mm of node,minimum width=55mm] (neuraltrain)
    {Differentiable ODE integration};
\node[flowbox,below=16mm of sindy] (equations)
    {Explicit equations\\$\dot y=\Theta(y)\widehat\Xi$};
\draw[flow] (node)--(neuraltrain);
\draw[flow] (sindy)--(equations);

\node[flowbox,below=17mm of neuraltrain,xshift=30mm,minimum width=75mm] (evaluation)
    {Reconstruction and future prediction\\relative errors and training time};
\draw[flow] (neuraltrain.south)--([xshift=-20mm]evaluation.north);
\draw[flow] (equations.south)--([xshift=20mm]evaluation.north);
\end{tikzpicture}%
}
\caption{\rev{Conceptual workflow. A NODE uses a neural network, parameterized by either an MLP or a KAN, within a differentiable ODE-training procedure. SINDy identifies an explicit sparse system.}}
\label{fig:learning-workflow-ode}
\end{figure}

\paragraph{\rev{Neural ODEs.}}
\rev{A Neural ODE (NODE) represents the unknown vector field by a neural network with trainable parameters $\theta$,}
\[
\frac{d y}{d t} \approx \hat{f}_{\theta}(y) = \mathrm{NN}(y; \theta).
\] 
\rev{Starting from $y(0)=y_0$, a differentiable ODE solver generates the trajectory $y_\theta(t)$, and the parameters are learned by minimizing \eqref{eq:loss}. We use a multilayer perceptron (MLP) as the standard neural parameterization.}

\paragraph{\rev{KANODEs.}}
\rev{A KANODE replaces the MLP in the same Neural ODE framework with a Kolmogorov--Arnold Network (KAN), whose learnable univariate spline functions parameterize the vector field }\rev{\cite{liu2024kan,koenig2024kan}}\rev{,}
\[
\frac{d y}{d t} \approx \hat{f}_{\theta}(y)
    = \mathrm{KAN}(y; \theta).
\]
\rev{The spline coefficients form $\theta$, while trajectory generation and training follow the same differentiable integration and loss used for the MLP-based NODE.}

\paragraph{\rev{SINDy.}}
\rev{Sparse Identification of Nonlinear Dynamics (SINDy) seeks an explicit governing equation from a prescribed library of candidate functions \cite{brunton2016discovering,de2020pysindy},}
\[
\frac{d y}{d t} \approx f_{\Xi}(y) = \Theta(y) \Xi,
\]
\rev{where $\Theta(y)$ contains candidate terms and the sparse coefficient matrix $\Xi$ is estimated by sparse regression, such as sequential thresholded least squares. Because numerical differentiation is sensitive to noise, regularization or trajectory smoothing may be applied before identification. 
}

\section{Learning Hidden Local Eco-Epidemic Coupling}

\rev{Finally, we present a  partial-learning framework for inferring hidden local eco--epidemic coupling.} We consider a spatio-temporal model by incorporating diffusion operators in one- and two-dimensional formulations. A finite difference method is applied to construct a discrete time-dependent system. 
\rev{Within this framework, the diffusion process is assumed to be known and is incorporated explicitly into the differentiable time-integration loop, while the hidden local eco--epidemic coupling is learned from solution data defined on the spatial grid. 
For efficient training on finer spatial grids, we formulate and train a reduced-order model based on Proper Orthogonal Decomposition (POD) and the QR-based variant of the Discrete Empirical Interpolation Method (QDEIM)~\cite{chaturantabut2010nonlinear, drmac2016new}. POD constructs a low-dimensional representation of the spatial state and provides low-rank filtering of noisy data, while DEIM hyper-reduces the nonlinear reaction by evaluating the learned coupling at a small set of reaction-informed nodes. The projected known diffusion operator is applied explicitly during reduced-space training, substantially reducing the cost of the forward and backward computations, while retaining accurate reconstruction on the original grid. We also use SINDy, introduced earlier as a main discovery component, synergistically after NODE training to obtain an optional explicit sparse representation of the learned reaction--diffusion system.}

\subsection{Spatio-Temporal Model with Finite Difference Diffusion}

To model diffusion in one or two spatial dimensions and discretize using the finite difference method. For a one-dimensional (1D) domain of length $X$ discretized into $N_x$ points with spacing $\Delta x = X/(N_x-1)$, the diffusion term
\[
\nabla (\kappa \nabla y) = \frac{\partial}{\partial x} \Big( \kappa \frac{\partial y}{\partial x} \Big)
\]
at interior points is approximated by the standard central difference
\[
\frac{\partial}{\partial x} \Big( \kappa \frac{\partial y}{\partial x} \Big) \Big|_{x_i} \approx \frac{\kappa_{i+1/2} (y_{i+1}-y_i) - \kappa_{i-1/2} (y_i - y_{i-1})}{\Delta x^2}, 
\quad i = 1, \dots, N_x-2,
\]
where $\kappa_{i \pm 1/2}$ denotes the diffusion coefficient at cell interfaces. 
\rev{The interface coefficients are computed using the harmonic mean, applied separately along each coordinate direction in two dimensions. For homogeneous diffusion, this reduces to the prescribed constant coefficient.} 
Homogeneous Neumann (zero-flux) boundary conditions are imposed on the domain boundary.

For a two-dimensional rectangular domain of size $L_x \times L_y$ discretized into $N_x \times N_y$ grid points
\[
\nabla (\kappa \nabla y) = 
\frac{\partial}{\partial x_1} \Big( \kappa \frac{\partial y}{\partial x_1} \Big) 
+ \frac{\partial}{\partial x_2} \Big( \kappa \frac{\partial y}{\partial x_2} \Big),
\quad x = (x_1, x_2), 
\]
the discrete Laplacian at interior nodes $(i,j)$ is approximated by
\[
\frac{\partial}{\partial x_1} \Big( \kappa \frac{\partial y}{\partial x_1} \Big) \Big|_{x_{ij}} 
\approx
\frac{\kappa_{i+1/2,j} (y_{i+1,j} - y_{i,j}) - \kappa_{i-1/2,j} (y_{i,j} - y_{i-1,j})}{\Delta x_1^2} ,
\]\[
\frac{\partial}{\partial x_2} \Big( \kappa \frac{\partial y}{\partial x_2} \Big) \Big|_{x_{ij}} 
\approx
\frac{\kappa_{i,j+1/2} (y_{i,j+1} - y_{i,j}) - \kappa_{i,j-1/2} (y_{i,j} - y_{i,j-1})}{\Delta x_2^2},
\]
where $(i,j)$ denotes the grid point in the 2D mesh with $x_{ij} = (x_{1,i}, x_{2,j})$, and $\kappa_{i\pm 1/2,j}$, $\kappa_{i,j\pm 1/2}$ are the diffusion coefficients at the cell interfaces in the $x_1$ and $x_2$ directions, respectively.

Including diffusion, the spatially discretized system \eqref{eq:mm} can be written as
\begin{equation}
\label{eq:spatio-temporal-fd}
\frac{d y}{d t} = Ly + f(y),
\end{equation}
where $y(t) \in \mathbb{R}^{M \times N}$ denotes the state vector collecting all $M$ components at $N$ spatial grid points ($N=N_x$ in 1D and $N=N_x\times N_y$ in 2D), $L$ is the discrete Laplacian obtained from a finite difference discretization, and $f(y)$ represents the local reaction term.
We note that in some cases, cross-diffusion terms can be incorporated to account for interactions between components.  In this case, we assume no cross-diffusion between components, so the resulting discrete Laplacian matrix $L$ has a block-diagonal structure, with each block corresponding to the diffusion of a single component. The resulting system of ODEs \eqref{eq:spatio-temporal-fd} can be integrated using the standard time-stepping schemes described above.
In this work, we consider systems with  small diffusion and employ explicit time-stepping schemes for temporal integration.

\subsection{\rev{Learning Hidden Local Coupling in a Reduced-Order Space}}

We now describe the learning strategy used to identify unknown local reaction coupling in the spatio-temporal model. Our primary focus is on learning the local reaction term, rather than fully resolving the coupled spatio-temporal behavior. Accordingly, we assume that the diffusion structure, or spatial connectivity, is known and fixed. Based on this assumption, we propose a partial learning framework that aims to identify hidden local interactions among the components of the system. This strategy enables us to infer the essential dynamics while leveraging prior knowledge of the spatial diffusion operator, thereby reducing the complexity of the learning problem and improving interpretability. Moreover, by concentrating on the local coupling terms, the proposed approach is able to capture subtle interactions that may be overlooked by fully data-driven models.

In this framework, the diffusion operator in \eqref{eq:spatio-temporal-fd} is explicitly incorporated into the model and treated as fixed, while the reaction term $f(y)$ is approximated by a neural network. This separation allows the learning process to focus exclusively on the unknown local nonlinear coupling between components, further reducing model complexity and enhancing interpretability. The resulting system is formulated as follows
\[
\frac{d y}{d t} = L y + \hat f_\theta(y),
\]
where $L$ denotes the known discrete Laplacian, and $\hat f_\theta$ is a neural network parameterized by $\theta$.

In the training process, the spatio-temporal solution data are normalized using a global mean and standard deviation computed over all time steps and spatial locations.
The training strategy explicitly separates the known diffusion dynamics from the unknown reaction processes, allowing the learning algorithm to focus on identifying the local nonlinear coupling between system components.
Since the diffusion term is treated deterministically and applied component-wise using a preassembled Laplacian matrix, the resulting NODEs preserve the underlying physical structure of the spatio-temporal model while benefiting from the flexibility of neural network-based closures.

\rev{Learning reaction--diffusion dynamics directly on a fine spatial grid is computationally expensive because every time-integration stage and gradient evaluation evolves all $N$ spatial nodes and evaluates the nonlinear reaction over the complete domain. This cost is particularly significant for two-dimensional grids and repeated trajectory rollouts during NODE training. To reduce both parts of the computation, we combine Proper Orthogonal Decomposition (POD) with the Discrete Empirical Interpolation Method (DEIM) and use its QR-based node-selection variant, QDEIM~\cite{chatterjee2000introduction,chaturantabut2010nonlinear,drmac2016new}. POD reduces the state and the known diffusion operator, whereas QDEIM hyper-reduces the evaluation of the learned nonlinear reaction.

Let, $Y^n\in\mathbb{R}^{N\times M}$ denote a normalized solution snapshot at time $t^n$, with $N$ spatial nodes and $M$ components. We select $s$ snapshots only from the training interval and concatenate their component fields into
\begin{equation}
X_y=[Y^{n_1},\ldots,Y^{n_s}]\in\mathbb{R}^{N\times Ms},
\qquad
X_y=U_y\Sigma_yV_y^{\top}.
\label{eq:pod-snapshot-svd}
\end{equation}

Let, $u_1^y,\ldots,u_r^y$ denote the leading $r$ left singular vectors of $X_y$. The shared spatial POD basis is defined as $\Phi_r=[u_1^y,\ldots,u_r^y]\in\mathbb{R}^{N\times r}$. The component fields are projected and reconstructed as $Z^n=\Phi_r^{\top}Y^n$ and $\widehat Y^n=\Phi_rZ^n$, respectively. Because $\Phi_r$ is constructed exclusively from training snapshots, the unseen future interval does not enter the reduced basis. Truncating the SVD also provides implicit low-rank filtering, where dominant coherent spatial patterns are retained, while low-energy modes that may contain measurement noise or small numerical fluctuations are removed. The rank therefore balances compression and denoising against the need to retain dynamically relevant spatial structure.

Galerkin projection of the full-order model onto this basis gives
\begin{equation}
\frac{dZ}{dt}
=L_rZ+\widehat f^{\mathrm{POD}}_{r,\theta}(Z),
\qquad
L_r=\Phi_r^{\top}L\Phi_r,
\qquad
\widehat f^{\mathrm{POD}}_{r,\theta}(Z)
=\Phi_r^{\top}\widehat f_\theta(\Phi_rZ),
\label{eq:pod-reduced-model}
\end{equation}
where $Z(t)\in\mathbb{R}^{r\times M}$ contains the reduced coefficients and $L_r\in\mathbb{R}^{r\times r}$ is the fixed projected diffusion operator. We note that \eqref{eq:pod-reduced-model} reduces the state and diffusion calculations, but evaluating $\widehat f_\theta(\Phi_rZ)$ still reconstructs the state and applies the reaction network at all $N$ nodes. QDEIM is introduced to remove this remaining full-grid nonlinear cost.

We first estimate reaction snapshots on the same training interval by subtracting the known diffusion contribution from finite-difference temporal derivatives, $R^n\approx\dot Y^n-LY^n$. After consistent normalization, the component fields are concatenated into a reaction snapshot matrix $X_f\in\mathbb{R}^{N\times Ms}$. Let $u_1^f,\ldots,u_q^f$ denote the leading $q$ left singular vectors of $X_f$. The reaction basis is defined as $U_f=[u_1^f,\ldots,u_q^f]\in\mathbb{R}^{N\times q}$. QDEIM applies column-pivoted QR factorization to $U_f^{\top}$ and selects $q$ informative, linearly independent spatial indices. Let $P=[e_{p_1},\ldots,e_{p_q}]\in\mathbb{R}^{N\times q}$ be the associated selection matrix. Enforcing interpolation at these nodes gives the DEIM approximation 
$f(Y)\approx U_f(P^{\top}U_f)^{-1}P^{\top}f(Y)$.

\begin{figure}[h!]
\centering
\resizebox{0.6\linewidth}{!}{%
\begin{tikzpicture}[
    font=\normalsize, >=Latex, node distance=13mm and 15mm,
    flowbox/.style={draw=black,text=black,rounded corners,align=center,
        minimum height=14mm,minimum width=40mm,fill=black!4},
    optional/.style={draw=black,text=black,dashed,rounded corners,align=center,
        minimum height=16mm,minimum width=48mm,fill=black!7},
    flow/.style={->,thick,draw=black}
]
\node[flowbox] (data) {1D/2D spatial trajectories\\normalization};
\node[optional,right=of data] (reduction) {Reduced space: POD--DEIM\\$z=\Phi_r^{\top}y$};
\node[flowbox,below=of reduction,xshift=-30mm,minimum width=52mm] (model)
    {Partial-learning\\known diffusion + learned reaction \\ NODE\\$\dot z=L_rz+\widehat f_{r,\theta}(z)$};
\draw[flow] (data)--(reduction);
\draw[flow] (reduction.south)--(model.north);
\draw[flow,dashed] (data.south) to[out=-55,in=155]
    node[left,text=black]{Full space ($\Phi_r=I$)} (model.west);

\node[flowbox,below=of model,minimum width=52mm] (trajectory)
    {Differentiable ODE integration\\ with learned reaction $\widehat f_{r,\theta}$\\fine-grid trajectory $\widehat y=\Phi_r z$};
\draw[flow] (model)--(trajectory);

\node[flowbox,below=15mm of trajectory,xshift=-34mm,minimum width=55mm] (prediction)
    {Prediction\\diffusion + learned reaction};
\node[optional,below=15mm of trajectory,xshift=34mm,minimum width=58mm] (sindy)
    {Optional SINDy identification\\fit the learned NODE trajectory\\explicit equations};
\draw[flow] (trajectory.south)--([xshift=-17mm]prediction.north);
\draw[flow] (trajectory.south)--([xshift=17mm]sindy.north);

\node[flowbox,below=17mm of prediction,xshift=34mm,minimum width=75mm] (evaluation)
    {Model assessment\\training and unseen-future errors\\training time};
\draw[flow] (prediction.south)--([xshift=-17mm]evaluation.north);
\draw[flow] (sindy.south)--([xshift=17mm]evaluation.north);
\end{tikzpicture}
}
\caption{\rev{Workflow for learning hidden local eco--epidemic coupling. The diffusion operator is prescribed, while differentiable ODE integration learns the local reaction. This trajectory is used either directly for space--time hybrid prediction or as input to an optional SINDy stage that identifies an explicit sparse model. Both outputs are assessed over the training interval and the unseen future horizon.}}
\label{fig:learning-workflow-hidden}
\end{figure}

For online training, we precompute the sampled state basis $S_r=P^{\top}\Phi_r\in\mathbb{R}^{q\times r}$ and the reduced interpolation operator $B_r=\Phi_r^{\top}U_f(P^{\top}U_f)^{-1}\in\mathbb{R}^{r\times q}$. At each ODE evaluation, $S_rZ$ forms the state only at the selected nodes, the same local reaction network is evaluated at those $q$ nodes, and $B_r$ maps its sampled output directly into the POD coordinates. The POD--QDEIM model is therefore
\begin{equation}
\frac{dZ}{dt}
=L_rZ+\widehat f^{\mathrm{QDEIM}}_{r,\theta}(Z),
\qquad
\widehat f^{\mathrm{QDEIM}}_{r,\theta}(Z)
=B_r\widehat f_\theta(S_rZ).
\label{eq:qdeim-reduced-model}
\end{equation}
Because the reaction basis is constructed after the known diffusion has been removed, the selected nodes target the spatial signatures of the hidden local coupling rather than the prescribed transport process. The training loop evolves $rM$ reduced variables, applies the known diffusion explicitly through $L_rZ$, and evaluates the neural reaction at $q$ nodes instead of all $N$ nodes. The matrices $L_r$, $S_r$, and $B_r$ remain fixed during optimization. This reduces the cost of both the differentiable rollout and backpropagation while preserving the known diffusion within the POD subspace and concentrating the neural model on the unknown reaction. We train the POD--QDEIM system using the same windowed trajectory loss as the full-space NODE. After training, the predicted fine-grid trajectory is reconstructed as $\widehat Y=\Phi_rZ$, and both training and full-trajectory errors are evaluated in full space. }

\rev{Figure~\ref{fig:learning-workflow-hidden} summarizes the complete procedure.
We note that  (i) the prescribed diffusion operator is separated from the unknown reaction and evaluated explicitly at every stage of the differentiable training loop; (ii) the POD basis is built only from training data; (iii) QDEIM is constructed from reaction snapshots after diffusion has been removed, so its interpolation nodes target the complex hidden local interaction rather than the already known diffusion process; and (iv) SINDy, introduced above as a principal component for direct discovery of average ODE dynamics, is optionally applied after NODE training to extract an explicit sparse representation from the learned smoother trajectory. 
The reduced-order formulation enables efficient learning on fine spatial grids. Instead of evolving the full state $Y\in\mathbb{R}^{N\times M}$ during every differentiable rollout, training evolves the reduced coefficients $Z\in\mathbb{R}^{r\times M}$. The known diffusion contribution is evaluated using the fixed reduced operator $L_r$, while the neural reaction is evaluated at only $q$ QDEIM nodes rather than at all $N$ spatial locations. This reduction decreases the computational cost of training. At the same time, the learned reduced trajectory is reconstructed as $\widehat Y=\Phi_rZ$, allowing prediction errors to be evaluated at the original spatial resolution. The reduced-order formulation accelerates learning while retaining the prescribed diffusion physics and the accuracy needed to resolve complex hidden local interactions.}

\rev{The present formulation assumes that the physical diffusion operator $L$ is known. If an approximate operator $\widetilde L$ is used instead, the discrepancy $(L-\widetilde L)Y$ enters the residual assigned to the unknown reaction. The neural model may partially compensate for this discrepancy, but the recovered reaction can become biased. This source of error is distinct from the reduced-space approximation introduced by $L_r=\Phi_r^{\top}L\Phi_r$, which trades POD projection error for faster training and low-rank denoising. Joint estimation of the diffusion coefficients and hidden reaction remains an important direction for future work.}

\section{Numerical Experiments}

We consider three datasets generated by numerically solving the ODE system~\eqref{eq:mm} using the right-hand sides $f_{\mathrm{SIR}}$ in~\eqref{sir}, $f_{\mathrm{LV}}$ in~\eqref{lv}, and the coupled model $f_{\mathrm{LVSIS}}$ in~\eqref{LVSIS}. 
The SIR, LV, and LVSIS systems are integrated over a fixed time interval with prescribed initial conditions to obtain the reference trajectories. \begin{itemize}
\item \textit{SIR: }
The \textit{SIR model} has state vector $y=(u_s,u_i,u_r)^\top$, with dynamics governed by \eqref{sir}, parameters $\sigma = 0.5$ and $\omega = 0.1$; and initial condition $y(0)=[0.99,0.01,0.0]^\top$.
\item \textit{LV: }
The \textit{Lotka--Volterra (predator--prey) model} has state vector $y=(u,v)^\top$, representing prey ($u$) and predator ($v$) populations. Its dynamics follow \eqref{lv} with parameters $\alpha = 0.5$, $\beta = 0.8$, $\delta = 0.4$, and $\gamma = 0.3$; initial condition $y(0)=[0.5,0.5]^\top$; and a small prey mortality $\eta = 0.05$.
\item \textit{LVSIS: }
The \textit{Lotka--Volterra--SIS (eco-epidemiological) model} has state vector $y=(u_s, u_i, v_s, v_i)^\top$ and dynamics given by \eqref{LVSIS}. The initial condition is $y(0)=[0.7, 0.3, 0.6, 0.4]^\top$. Prey birth rates are set to $\alpha_s = 0.04$ and $\alpha_i = 0.02$; and prey and predator death rates are $\gamma^u_s = \gamma^u_i = 0.02$ and $\gamma^v_s = \gamma^v_i = 0.01$, respectively. Coupling parameters for predator--prey interactions are $\beta_{s}^s = \beta_{s}^i = \beta_{i}^s = \beta_{i}^i = 0.02$, and for cross-species interactions $\delta_{s}^s = \delta_{s}^i = \delta_{i}^s = \delta_{i}^i = 0.04$. Transmission and recovery rates are $\sigma^s_{i} = \sigma^s_{vi} = \sigma^v_{si} = \sigma^{vs}_{vi} = 0.01$ and $\omega_i = \omega_{vi} = 0.01$. 
\end{itemize}

We begin simulations at $t=0$ and evolve the system up to $T_{\max}$, with $T_{\max} = 60, 180$, and $1971$ for the SIR, LV, and LVSIS datasets, corresponding to the right-hand sides $f_{\mathrm{SIR}}$ in~\eqref{sir}, $f_{\mathrm{LV}}$ in~\eqref{lv}, and $f_{\mathrm{LVSIS}}$ in~\eqref{LVSIS}, respectively. Each simulation is discretized into $N_t = 300$ time steps, producing a dataset that is subsequently split into training and testing subsets. We consider three training proportions $p = 1/3, 1/2$, and $2/3$, corresponding to the first $N = 100, 150$, and $200$ time steps. The resulting time series $\tilde{Y}$, consisting of $N$ points, serves as input for the data-driven discovery methods, including Neural ODEs, KANODEs, and SINDy. 
The ODE systems in \eqref{eq:mm} are numerically solved using the SciPy library \cite{virtanen2020scipy}.

We consider the performance of three techniques to learn ecological and epidemic dynamics based on the given dataset. 
For the KAN architecture, the number of trainable parameters is 360 for the SIR dataset, 
160 for the LV dataset, and 640 for the LVSIS dataset ($N_{comp}$ - \rev{2 $N_{comp}$} - $N_{comp}$).
For NODE$_1$, the parameter counts are 4,611 for the SIR, 4,482 for the LV, and 4,740 for the LVSIS ($N_{comp}$-64-$N_{comp}$). 
For NODE$_2$, the models have 33,539 parameters for SIR, 33,410 for LV, and 33,668 for LVSIS, 
with a deeper architecture ($N_{comp}$-64-128-128-64-$N_{comp}$).

\rev{All machine-learning models are trained for $500$ epochs using the Adam optimizer. For NODE$_1$ and NODE$_2$, the activation function after each hidden layer is $\tanh$. Training is performed using mini-batches of short trajectory windows. At each optimization step, the starting indices of the windows are sampled randomly from the training interval. 
We use a learning rate of $10^{-3}$, a batch size of $100$, a window length of $10$ time steps, and $10$ optimization steps per epoch. For the results with diffusion, we use a learning rate of $5.0\times 10^{-4}$, a batch size of $20$, a window length of $10$ time steps, and $10$ optimization steps per epoch \cite{lee2022structure}. 
Following the KAN formulations~\cite{liu2024kan,koenig2024kan}, each KAN layer uses a fixed grid with $5$ intervals on the initial range $[-1,1]$ and cubic B-splines of order $3$. The base branch uses the SiLU activation, the base and spline scales are both initialized to $1$, the spline initialization-noise scale is $0.1$, and a separately trainable spline-scaling coefficient is used for every connection. The KAN grid is not updated during training. 
For direct SINDy identification of the ODE systems, temporal derivatives are computed using second-order finite differences, and the candidate library contains polynomial terms up to degree $2$. Sparse regression is performed using Sequential Thresholded Least Squares (STLSQ) with threshold $10^{-4}$ and ridge regularization parameter $10^{-3}$. 
We use the PySINDy library \cite{de2020pysindy} for sparse system identification and PyTorch \cite{paszke2019pytorch} for training the neural network models.
All experiments are performed on a MacBook Air equipped with a 10-core Apple M4 processor and 16~GB of memory.}

The prediction accuracy is quantified using relative $\ell^2$ errors expressed as percentages. 
At selected time indices $i \in \mathcal{I}$, the pointwise error is computed as
\[
e^n
= \frac{\lVert y_{\mathrm{pred}}^n - y^n \rVert_2}
{\lVert y^n \rVert_2} \times 100 \%,
\]
In addition, we report relative $\ell^2$ errors on the training subset and on the full dataset
\[
e_{\text{Train}} = \frac{
\left( \sum_{n=1}^{N}
\left\lVert
y_{\theta}^n - y^n
\right\rVert^2
\right)^{1/2}
}{
\left( \sum_{n=1}^{N}
\left\lVert
y^n
\right\rVert^2
\right)^{1/2}
}
\times 100 \%, 
\quad 
e_{\text{Full}} = \frac{
\left( \sum_{n=1}^{N_t}
\left\lVert
y_{\theta}^n - y^n
\right\rVert^2
\right)^{1/2}
}{
\left( \sum_{n=1}^{N_t}
\left\lVert
y^n
\right\rVert^2
\right)^{1/2}
}
\times 100 \%. 
\]

\subsection{Learning Ecological and Epidemic Dynamics}

\rev{In this section, we evaluate how accurately the temporal learning methods recover and predict the dynamics of the three benchmark systems. We assess temporal prediction using the training and full-horizon errors in Tables~\ref{tab:temporal-sir}--\ref{tab:temporal-lvsis}. The noise-free and noisy-data trajectory plots, together with the detailed SINDy equations, are provided in Appendix~\ref{app:temporal-supplement}.
}

\begin{table}[h!]
\centering
\small
\begin{tabular}{|l|cc|cc|cc|}
\hline
\multicolumn{7}{|c|}{SIR} \\
\hline
 & \multicolumn{2}{c|}{1/3} & \multicolumn{2}{c|}{1/2} & \multicolumn{2}{c|}{2/3} \\
 & Train & Full & Train & Full & Train & Full \\
\hline
\multicolumn{7}{|c|}{Original data} \\
\hline
SINDy & 0.06 & 1.50 & 0.07 & 0.54 & 0.07 & 0.13 \\
DMD & 13.32 & 679.0 & 11.13 & 62.20 & 13.17 & 20.89 \\
EDMD & 6.02 & 28751.2 & 4.57 & 242.4 & 5.76 & 21.66 \\
NODE$_1$ & 2.72 $\pm$ 1.25 & 677.6 $\pm$ 44.79 & 1.45 $\pm$ 0.87 & 7.17 $\pm$ 1.66 & 1.16 $\pm$ 1.10 & 1.28 $\pm$ 0.93 \\
NODE$_2$ & 2.10 $\pm$ 0.91 & 412.0 $\pm$ 230.3 & 1.52 $\pm$ 0.81 & 3.09 $\pm$ 0.66 & 0.85 $\pm$ 0.98 & 0.80 $\pm$ 0.60 \\
KANODE & 0.09 $\pm$ 0.05 & 89.08 $\pm$ 43.14 & 0.06 $\pm$ 0.06 & 20.63 $\pm$ 6.76 & 0.46 $\pm$ 0.37 & 1.95 $\pm$ 0.26 \\
\hline
\multicolumn{7}{|c|}{Noisy data ($\sigma=0.01$)} \\
\hline
SINDy & - & - & - & - & - & - \\
DMD & 17.72 & 106.8 & 13.66 & 48.85 & 12.24 & 18.35 \\
EDMD & 8.61 & 51.27 & 7.94 & 42.93 & 7.34 & 12.07 \\
NODE$_1$ & 10.41 $\pm$ 0.40 & 32.69 $\pm$ 3.45 & 8.15 $\pm$ 1.26 & 6.19 $\pm$ 0.85 & 5.96 $\pm$ 1.81 & 4.44 $\pm$ 0.96 \\
NODE$_2$ & 8.67 $\pm$ 4.74 & 39.53 $\pm$ 2.01 & 10.59 $\pm$ 2.12 & 7.93 $\pm$ 1.59 & 7.48 $\pm$ 0.92 & 5.32 $\pm$ 0.63 \\
KANODE & 11.01 $\pm$ 0.94 & 41.12 $\pm$ 3.65 & 7.94 $\pm$ 1.46 & 39.89 $\pm$ 27.84 & 8.04 $\pm$ 1.11 & 7.18 $\pm$ 0.57 \\
\hline
\multicolumn{7}{|c|}{Noisy data ($\sigma=0.05$)} \\
\hline
SINDy & - & - & - & - & - & - \\
DMD & 22.75 & 51.95 & 31.07 & 38.68 & 31.85 & 32.55 \\
EDMD & 23.55 & 51.28 & 27.41 & 36.23 & 27.30 & 29.04 \\
NODE$_1$ & 27.74 $\pm$ 0.43 & 42.76 $\pm$ 2.23 & 24.21 $\pm$ 0.97 & 24.96 $\pm$ 0.29 & 22.75 $\pm$ 1.19 & 17.76 $\pm$ 1.17 \\
NODE$_2$ & 23.86 $\pm$ 3.45 & 40.24 $\pm$ 1.95 & 18.41 $\pm$ 3.71 & 21.45 $\pm$ 0.74 & 24.38 $\pm$ 5.17 & 19.29 $\pm$ 3.35 \\
KANODE & 23.78 $\pm$ 3.75 & 42.34 $\pm$ 1.00 & 18.54 $\pm$ 5.42 & 21.40 $\pm$ 2.64 & 9.19 $\pm$ 3.58 & 9.07 $\pm$ 2.71 \\
\hline
\multicolumn{7}{|c|}{Noisy data ($\sigma=0.1$)} \\
\hline
SINDy & - & - & - & - & - & - \\
DMD & 36.43 & 68.05 & 44.63 & 54.64 & 47.25 & 48.70 \\
EDMD & 36.35 & 67.76 & 42.91 & 53.53 & 43.45 & 45.59 \\
NODE$_1$ & 34.35 $\pm$ 1.35 & 53.29 $\pm$ 0.15 & 29.12 $\pm$ 3.07 & 31.83 $\pm$ 1.52 & 27.24 $\pm$ 0.78 & 23.69 $\pm$ 0.56 \\
NODE$_2$ & 26.89 $\pm$ 12.90 & 52.23 $\pm$ 4.64 & 26.60 $\pm$ 11.17 & 30.06 $\pm$ 6.49 & 28.83 $\pm$ 6.59 & 25.42 $\pm$ 3.82 \\
KANODE & 41.36 $\pm$ 7.39 & 47.97 $\pm$ 0.39 & 35.96 $\pm$ 2.45 & 33.45 $\pm$ 0.48 & 28.35 $\pm$ 6.10 & 28.86 $\pm$ 8.04 \\
\hline 
\end{tabular}
\caption{\rev{Training-interval and full-trajectory relative errors (\%) for the SIR temporal system with noise-free data and additive Gaussian noise at $\sigma=0.01$, $0.05$, and $0.1$. Neural-model entries report mean $\pm$ standard deviation over random initializations; DMD and EDMD are deterministic fits, and SINDy is reported as a deterministic fit for noise-free data only.}}
\label{tab:temporal-sir}
\end{table}

\begin{table}[h!]
\centering
\small
\begin{tabular}{|l|cc|cc|cc|}
\hline
\multicolumn{7}{|c|}{LV} \\
\hline
 & \multicolumn{2}{c|}{1/3} & \multicolumn{2}{c|}{1/2} & \multicolumn{2}{c|}{2/3} \\
 & Train & Full & Train & Full & Train & Full \\
\hline
\multicolumn{7}{|c|}{Original data} \\
\hline
SINDy & 1.23 & 1.32 & 1.38 & 1.30 & 1.38 & 1.29 \\
DMD & 2.53 & 1.53 & 1.73 & 1.30 & 1.51 & 1.32 \\
EDMD & 1.16 & 0.63 & 0.65 & 0.56 & 0.59 & 0.60 \\
NODE$_1$ & 0.29 $\pm$ 0.07 & 0.28 $\pm$ 0.04 & 0.30 $\pm$ 0.06 & 0.33 $\pm$ 0.05 & 0.18 $\pm$ 0.01 & 0.25 $\pm$ 0.03 \\
NODE$_2$ & 0.29 $\pm$ 0.10 & 0.36 $\pm$ 0.15 & 0.26 $\pm$ 0.13 & 0.31 $\pm$ 0.13 & 0.32 $\pm$ 0.15 & 0.35 $\pm$ 0.15 \\
KANODE & 0.33 $\pm$ 0.03 & 1.84 $\pm$ 0.43 & 0.24 $\pm$ 0.07 & 0.90 $\pm$ 0.23 & 0.42 $\pm$ 0.12 & 0.37 $\pm$ 0.12 \\
\hline
\multicolumn{7}{|c|}{Noisy data ($\sigma=0.01$)} \\
\hline
SINDy & - & - & - & - & - & - \\
DMD & 5.72 & 3.72 & 6.53 & 4.52 & 6.74 & 5.17 \\
EDMD & 5.16 & 3.67 & 6.16 & 4.36 & 6.49 & 5.01 \\
NODE$_1$ & 1.02 $\pm$ 0.33 & 0.89 $\pm$ 0.36 & 1.19 $\pm$ 0.40 & 0.96 $\pm$ 0.42 & 0.96 $\pm$ 0.29 & 0.94 $\pm$ 0.27 \\
NODE$_2$ & 1.55 $\pm$ 0.23 & 4.24 $\pm$ 0.95 & 1.51 $\pm$ 0.30 & 1.65 $\pm$ 0.55 & 1.76 $\pm$ 0.63 & 1.62 $\pm$ 0.57 \\
KANODE & 0.90 $\pm$ 0.06 & 4.46 $\pm$ 0.76 & 1.06 $\pm$ 0.18 & 1.95 $\pm$ 0.83 & 1.19 $\pm$ 0.18 & 1.48 $\pm$ 0.14 \\
\hline
\multicolumn{7}{|c|}{Noisy data ($\sigma=0.05$)} \\
\hline
SINDy & - & - & - & - & - & - \\
DMD & 15.62 & 8.07 & 12.99 & 7.98 & 10.79 & 7.90 \\
EDMD & 15.76 & 8.13 & 13.07 & 8.02 & 10.81 & 7.91 \\
NODE$_1$ & 4.38 $\pm$ 0.23 & 8.12 $\pm$ 0.19 & 5.23 $\pm$ 0.74 & 3.75 $\pm$ 0.29 & 7.77 $\pm$ 0.13 & 5.88 $\pm$ 0.09 \\
NODE$_2$ & 8.76 $\pm$ 1.68 & 11.63 $\pm$ 1.45 & 5.56 $\pm$ 0.44 & 5.72 $\pm$ 0.79 & 7.72 $\pm$ 0.34 & 5.86 $\pm$ 0.22 \\
KANODE & 5.39 $\pm$ 0.61 & 10.20 $\pm$ 0.66 & 5.82 $\pm$ 0.65 & 6.50 $\pm$ 1.22 & 7.66 $\pm$ 0.73 & 5.92 $\pm$ 0.57 \\
\hline
\multicolumn{7}{c}{Noisy data ($\sigma=0.1$)} \\
\hline
SINDy & - & - & - & - & - & - \\
DMD & 16.14 & 8.25 & 13.23 & 8.09 & 10.95 & 8.00 \\
EDMD & 16.20 & 8.24 & 13.27 & 8.11 & 10.98 & 8.01 \\
NODE$_1$ & 9.50 $\pm$ 1.13 & 7.92 $\pm$ 1.48 & 11.43 $\pm$ 1.79 & 9.29 $\pm$ 0.82 & 9.86 $\pm$ 0.34 & 7.65 $\pm$ 0.65 \\
NODE$_2$ & 10.32 $\pm$ 3.64 & 11.51 $\pm$ 2.33 & 11.68 $\pm$ 3.74 & 10.35 $\pm$ 2.08 & 9.19 $\pm$ 1.53 & 8.80 $\pm$ 1.70 \\
KANODE & 10.10 $\pm$ 4.28 & 11.21 $\pm$ 0.76 & 9.18 $\pm$ 4.09 & 9.09 $\pm$ 2.81 & 7.23 $\pm$ 0.57 & 7.25 $\pm$ 0.73 \\
\hline
\end{tabular}
\caption{\rev{Training-interval and full-trajectory relative errors (\%) for the LV temporal system with noise-free data and additive Gaussian noise at $\sigma=0.01$, $0.05$, and $0.1$. Neural-model entries report mean $\pm$ standard deviation over random initializations; DMD and EDMD are deterministic fits, and SINDy is reported as a deterministic fit for noise-free data only.}}
\label{tab:temporal-lv}
\end{table}

\begin{table}[h!]
\centering
\small
\begin{tabular}{|l|cc|cc|cc|}
\hline
\multicolumn{7}{|c|}{LVSIS, noise-free} \\
\hline
 & \multicolumn{2}{c|}{1/3} & \multicolumn{2}{c|}{1/2} & \multicolumn{2}{c|}{2/3} \\
 & Train & Full & Train & Full & Train & Full \\
\hline
\multicolumn{7}{|c|}{Original data} \\
\hline
SINDy & 0.54 & 1.09 & 0.72 & 0.73 & 0.72 & 0.69 \\
DMD & 28.82 & 21.42 & 26.70 & 19.61 & 23.70 & 18.88 \\
EDMD & 7.56 & 56.87 & 6.29 & 6.80 & 6.04 & 5.96 \\
NODE$_1$ & 2.04 $\pm$ 1.02 & 7.03 $\pm$ 2.99 & 2.96 $\pm$ 1.24 & 2.80 $\pm$ 0.73 & 1.81 $\pm$ 0.16 & 1.74 $\pm$ 0.15 \\
NODE$_2$ & 4.33 $\pm$ 1.39 & 11.01 $\pm$ 0.85 & 3.50 $\pm$ 0.91 & 3.05 $\pm$ 1.42 & 2.89 $\pm$ 1.98 & 2.39 $\pm$ 1.42 \\
KANODE & 0.18 $\pm$ 0.06 & 11.96 $\pm$ 1.26 & 0.33 $\pm$ 0.18 & 1.50 $\pm$ 0.42 & 0.43 $\pm$ 0.25 & 0.54 $\pm$ 0.24 \\
\hline
\multicolumn{7}{|c|}{Noisy data ($\sigma=0.01$)} \\
\hline
SINDy & - & - & - & - & - & - \\
DMD & 26.37 & 18.58 & 23.80 & 16.72 & 20.07 & 15.57 \\
EDMD & 7.13 & 6.99 & 7.86 & 6.61 & 8.58 & 7.55 \\
NODE$_1$ & 5.67 $\pm$ 1.89 & 10.17 $\pm$ 0.13 & 4.68 $\pm$ 1.86 & 4.53 $\pm$ 1.74 & 5.25 $\pm$ 0.89 & 4.64 $\pm$ 0.65 \\
NODE$_2$ & 4.31 $\pm$ 1.63 & 8.80 $\pm$ 1.64 & 6.47 $\pm$ 3.70 & 5.88 $\pm$ 2.82 & 5.82 $\pm$ 1.83 & 5.06 $\pm$ 1.55 \\
KANODE & 2.41 $\pm$ 0.47 & 13.13 $\pm$ 2.71 & 1.71 $\pm$ 0.08 & 4.57 $\pm$ 1.03 & 4.00 $\pm$ 0.71 & 4.51 $\pm$ 0.67 \\
\hline
\multicolumn{7}{|c|}{Noisy data ($\sigma=0.05$)} \\
\hline
SINDy & - & - & - & - & - & - \\
DMD & 30.73 & 17.82 & 29.88 & 18.98 & 26.73 & 19.81 \\
EDMD & 26.75 & 16.05 & 26.53 & 17.06 & 23.68 & 17.63 \\
NODE$_1$ & 8.39 $\pm$ 1.96 & 19.82 $\pm$ 8.11 & 9.81 $\pm$ 3.04 & 13.16 $\pm$ 4.14 & 13.31 $\pm$ 4.92 & 10.87 $\pm$ 3.14 \\
NODE$_2$ & 17.59 $\pm$ 9.82 & 27.19 $\pm$ 11.02 & 15.52 $\pm$ 4.07 & 16.94 $\pm$ 5.46 & 13.62 $\pm$ 1.05 & 11.12 $\pm$ 0.76 \\
KANODE & 10.39 $\pm$ 3.03 & 27.58 $\pm$ 5.29 & 12.60 $\pm$ 3.82 & 18.31 $\pm$ 2.28 & 16.33 $\pm$ 3.75 & 16.83 $\pm$ 3.29 \\
\hline
\multicolumn{7}{|c|}{Noisy data ($\sigma=0.1$)} \\
\hline
SINDy & - & - & - & - & - & - \\
DMD & 43.36 & 23.37 & 36.48 & 22.92 & 30.63 & 22.83 \\
EDMD & 38.81 & 21.06 & 33.69 & 21.11 & 28.72 & 21.30 \\
NODE$_1$ & 19.95 $\pm$ 11.94 & 33.33 $\pm$ 2.91 & 22.84 $\pm$ 1.50 & 17.27 $\pm$ 1.23 & 19.82 $\pm$ 1.63 & 15.64 $\pm$ 1.16 \\
NODE$_2$ & 13.64 $\pm$ 3.19 & 30.85 $\pm$ 3.01 & 16.74 $\pm$ 4.94 & 20.17 $\pm$ 2.49 & 19.84 $\pm$ 0.58 & 16.63 $\pm$ 2.30 \\
KANODE & 33.69 $\pm$ 7.27 & 41.23 $\pm$ 1.25 & 16.36 $\pm$ 7.47 & 18.57 $\pm$ 4.17 & 13.41 $\pm$ 2.34 & 13.79 $\pm$ 2.39 \\
\hline
\end{tabular}
\caption{\rev{Training-interval and full-trajectory relative errors (\%) for the LVSIS temporal system with noise-free data and additive Gaussian noise at $\sigma=0.01$, $0.05$, and $0.1$. Neural-model entries report mean $\pm$ standard deviation over random initializations; DMD and EDMD are deterministic fits, and SINDy is reported as a deterministic fit for noise-free data only.}}
\label{tab:temporal-lvsis}
\end{table}

\rev{In addition to SINDy, NODEs, and KANODEs, the tables include Dynamic Mode Decomposition (DMD) and Extended Dynamic Mode Decomposition (EDMD) as standard non-neural baselines~\cite{schmid2010dynamic,williams2015koopman}. To state the two operator-based fits precisely, let $\bar y_n\in\mathbb{R}^{d}$ denote the normalized state at time $t_n$. We use the degree-one and degree-two observable maps
\begin{equation}
\psi_1(\bar y)
=
\begin{bmatrix}
1 & \bar y_1 & \cdots & \bar y_d
\end{bmatrix}^{\mathsf T},
\qquad
\psi_2(\bar y)
=
\begin{bmatrix}
1 & \bar y_1 & \cdots & \bar y_d
& \{\bar y_i\bar y_j\}_{1\le i\le j\le d}
\end{bmatrix}^{\mathsf T}.
\end{equation}
Thus, $\psi_2$ contains the constant, every normalized state component, and every distinct quadratic monomial. For either $q=1$ or $q=2$, define consecutive lifted snapshot matrices
\begin{equation}
Z_q
=
\begin{bmatrix}
\psi_q(\bar y_0) & \cdots & \psi_q(\bar y_{m-1})
\end{bmatrix},
\qquad
Z_q^{\mathrm{next}}
=
\begin{bmatrix}
\psi_q(\bar y_1) & \cdots & \psi_q(\bar y_m)
\end{bmatrix}.
\end{equation}
The full-rank least-squares evolution operator and the decoder to the normalized state are
\begin{equation}
K_q
=
\underset{K}{\operatorname{argmin}}\,
\left\lVert Z_q^{\mathrm{next}}-KZ_q\right\rVert_F^2
=
Z_q^{\mathrm{next}}Z_q^\dagger,
\qquad
C_q
=
\bar Y\left(Z_q^{\mathrm{all}}\right)^\dagger,
\end{equation}
\begin{equation}
\bar Y
=
\begin{bmatrix}
\bar y_0 & \cdots & \bar y_m
\end{bmatrix},
\qquad
Z_q^{\mathrm{all}}
=
\begin{bmatrix}
\psi_q(\bar y_0) & \cdots & \psi_q(\bar y_m)
\end{bmatrix},
\end{equation}
where ${}^\dagger$ denotes the Moore--Penrose pseudoinverse computed by singular value decomposition without prescribed rank truncation. Starting from $\hat z_0=\psi_q(\bar y_0)$, prediction proceeds recursively as
\begin{equation}
\hat z_{n+1}=K_q\hat z_n,
\qquad
\hat{\bar y}_n=C_q\hat z_n,
\end{equation}
followed by inversion of the componentwise normalization. The reported DMD baseline uses $q=1$; because $\psi_1$ includes a constant observable, it is an affine degree-one evolution model in the original normalized coordinates. The EDMD baseline uses $q=2$, so the lifted dynamics are linear while the observables are quadratic in the state. Consequently, DMD tests whether affine evolution in the original coordinates is sufficient, whereas EDMD tests whether a fixed quadratic lifting can represent the nonlinear dynamics.}

\rev{We consider the original noise-free trajectories together with noisy versions of these data. The noisy datasets are generated by adding independent zero-mean Gaussian noise to each reference trajectory}
\begin{equation}
\tilde{Y} = Y + \mathcal{N}(0, \sigma^2), \quad 
Y = \{ y(t_i) \}_{i=1}^N,
\end{equation}
\rev{where $Y$ denotes the original noise-free trajectory and $\sigma$ is the standard deviation of the added noise. We consider $\sigma\in\{0.01,0.05,0.1\}$. Results for the original noise-free data and the three noisy datasets are reported in Tables~\ref{tab:temporal-sir}--\ref{tab:temporal-lvsis}. Direct SINDy identification is reported only for noise-free data because finite-difference differentiation of the noisy trajectories produces spurious library terms and unstable coefficients under the fixed identification settings. The neural results are therefore compared with DMD and EDMD at every noise level. For NODE and KANODE models, the tables report the mean and standard deviation over three random initializations.}

\rev{Tables~\ref{tab:temporal-sir}--\ref{tab:temporal-lvsis} show that the appropriate baseline depends strongly on the dynamics and the prediction horizon. For noise-free LV data, the quadratic EDMD lifting is highly effective, with full-horizon errors of $0.60$--$0.63\%$, while DMD is also competitive; this confirms that standard operator-based techniques can be sufficient for regular oscillatory trajectories. Their limitations are clearer for SIR and LVSIS, where short SIR training windows lead to very large extrapolation errors for both baselines, and DMD gives full-horizon errors of approximately $19$--$21\%$ for noise-free LVSIS. The learned nonlinear models generally improve as the observed fraction increases, although no architecture dominates every case.}

\rev{The larger-noise results provide a more demanding robustness assessment. As $\sigma$ increases from $0.01$ to $0.1$, errors generally rise for every method, but the rate of degradation is system dependent. For LV with $p=2/3$, for example, the full-horizon errors at $\sigma=0.1$ are $7.65\%$, $8.80\%$, and $7.25\%$ for NODE$_1$, NODE$_2$, and KANODE, compared with about $8\%$ for DMD and EDMD. For LVSIS at the same noise level and training fraction, KANODE gives the smallest full-horizon error ($13.79\%$), followed by NODE$_1$ ($15.64\%$), whereas DMD and EDMD exceed $21\%$. The standard deviations also widen in several high-noise cases, especially for shorter training windows, showing that robustness depends on both noise amplitude and initialization. }
\rev{A detailed wall-clock comparison of SINDy, NODE$_1$, NODE$_2$, and KANODE, together with an LVSIS accuracy--time checkpoint study, is provided in Appendix~\ref{app:temporal-timing}.}

\subsection{Learning Hidden Local Coupling in Spatio-Temporal Systems}

\begin{figure}[h!]
\centering
\includegraphics[width=0.45\textwidth]{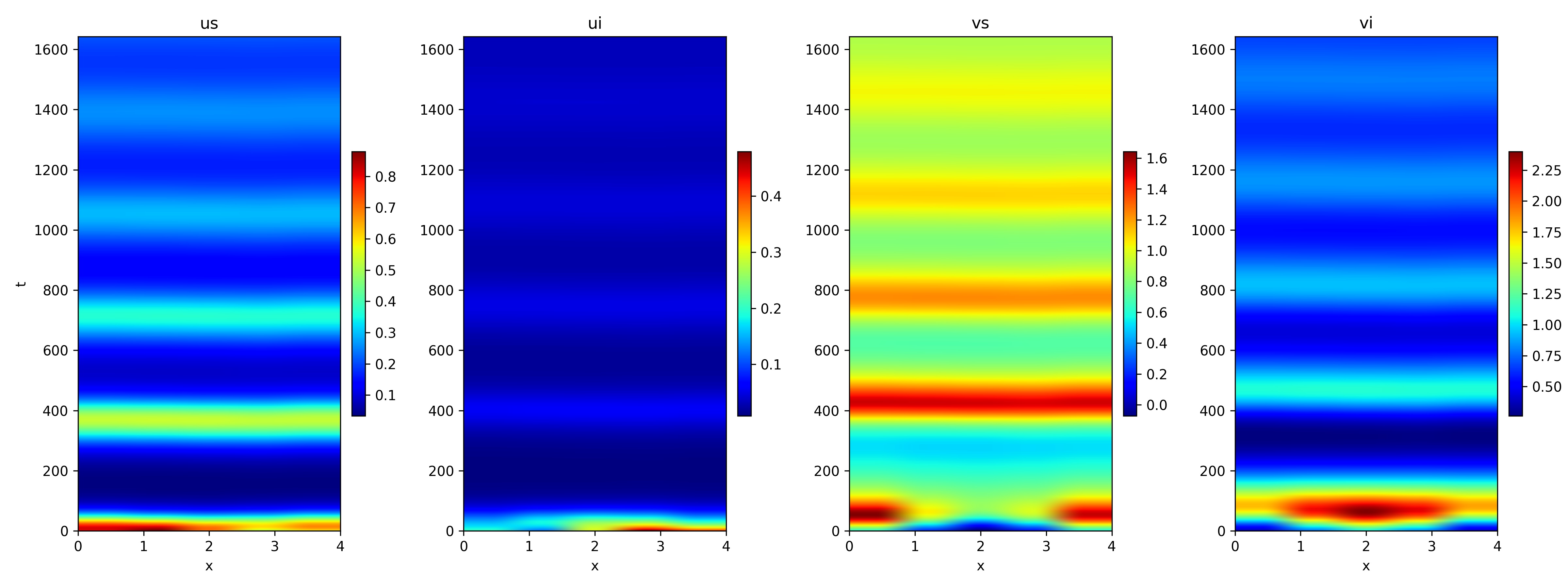} \ \ \ \ 
\includegraphics[width=0.45\textwidth]{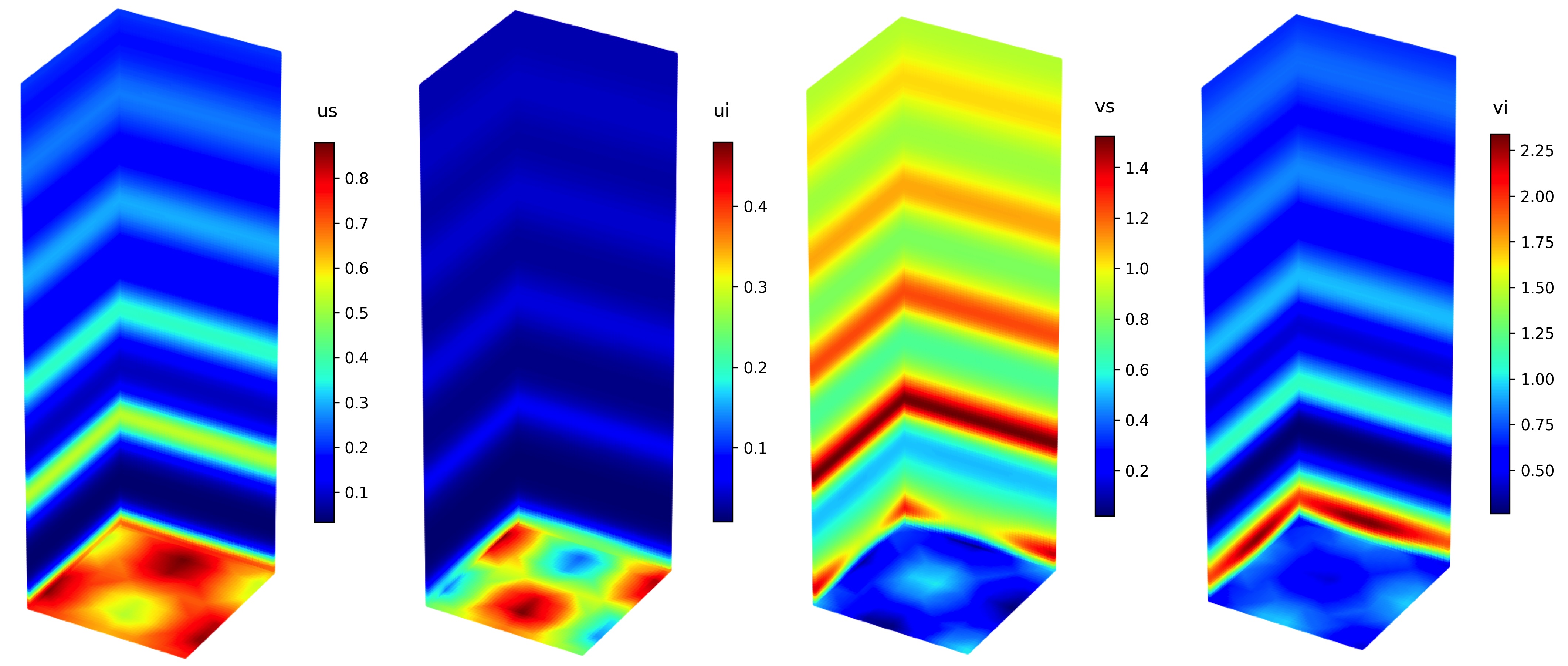}
\caption{\rev{LVSIS mode.  Reference solution. Spatio-temporal model (left: 1D, right:2D).}}
\label{fig:test-lvsis-ref}
\end{figure}

Finally, we consider the spatio-temporal dynamics of the LVSIS model. Figure~\ref{fig:test-lvsis-ref} presents the reference solutions in one and two spatial dimensions (1D and 2D). Our objective is to learn the local coupling induced by the LVSIS reaction terms, assuming the diffusion process is known. The diffusion is taken to be homogeneous, with diffusion coefficient $\kappa = 10^{-5}$ for each component, and is discretized using a finite difference method on uniform grids consisting of 10 nodes in the 1D case and $7 \times 7$ nodes in the 2D case. We consider the semi-discrete system \eqref{eq:spatio-temporal-fd} on the spatial domain $[0,10]^d$, for $d=1,2$. The parameters in the reaction term $f_{\mathrm{LVSIS}}$ are chosen to be consistent with those used in the previous experiments.

\begin{figure}[h!]
\centering
\setlength{\tabcolsep}{4pt}
\renewcommand{\arraystretch}{1.2}
\begin{tabular}{c c c c c}
 & 1/3 & 1/2 & 2/3 \\
\\
\centering\rotatebox{90}{Dynamics} & 
\includegraphics[width=0.3\textwidth]{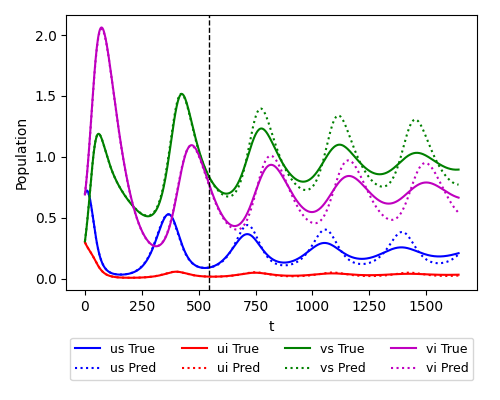} &
\includegraphics[width=0.3\textwidth]{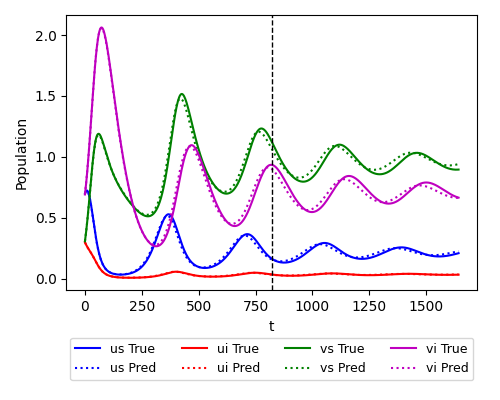} &
\includegraphics[width=0.3\textwidth]{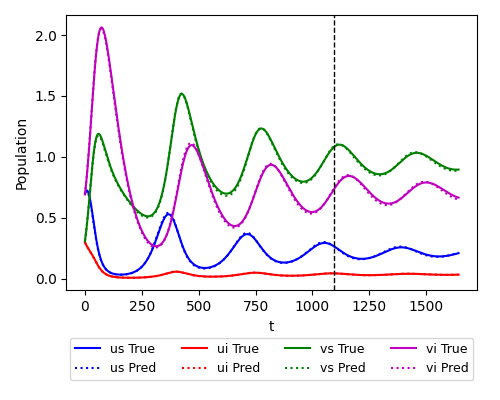} 
\\
\centering\rotatebox{90}{Prediction} & 
\includegraphics[width=0.3\textwidth]{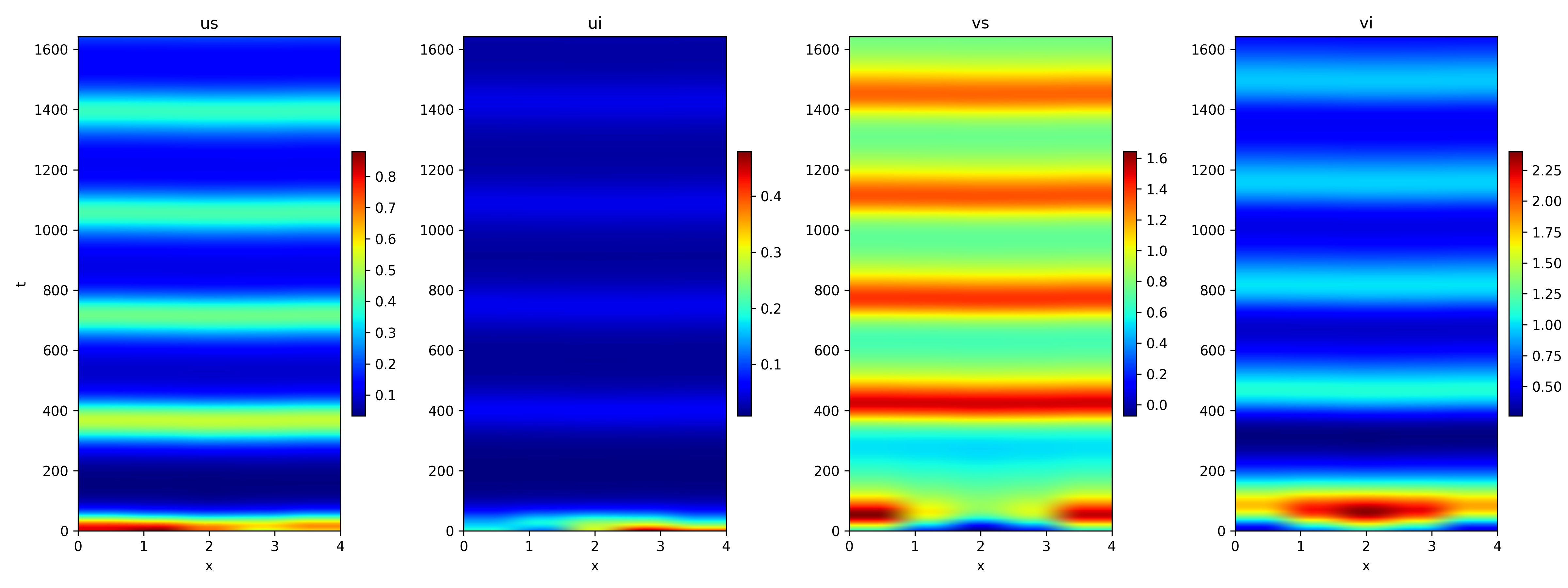} &
\includegraphics[width=0.3\textwidth]{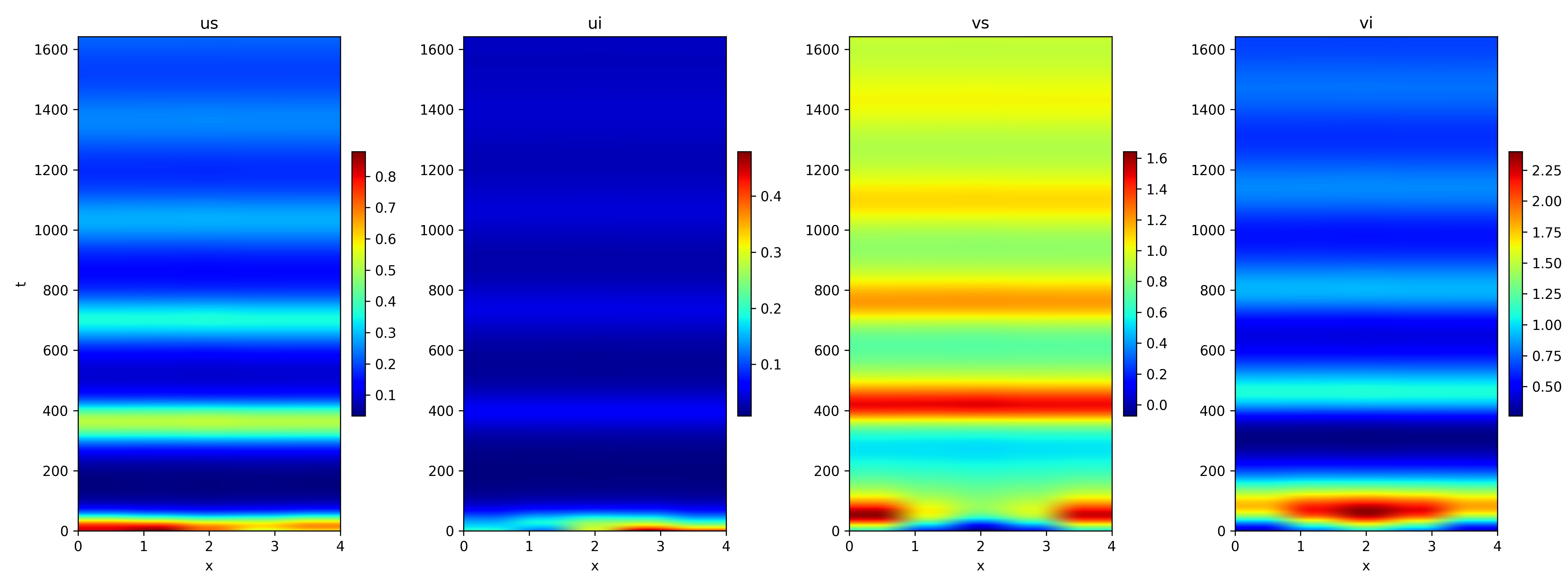} &
\includegraphics[width=0.3\textwidth]{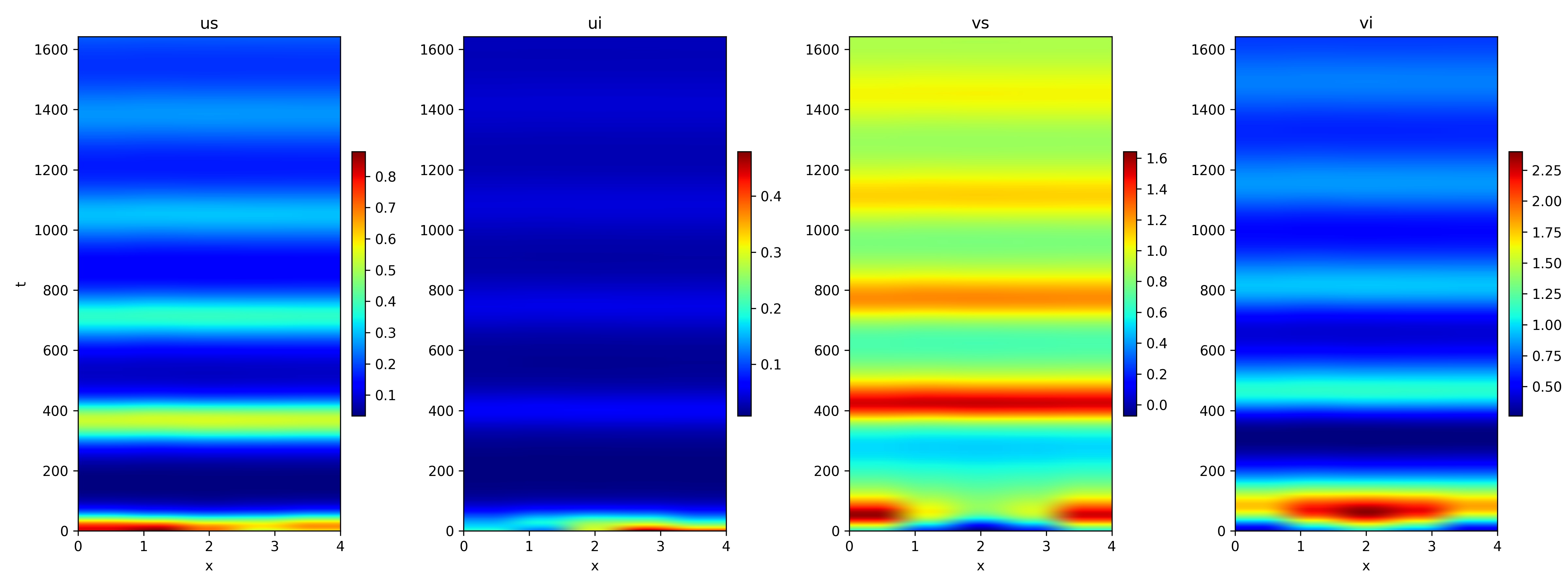} 
\end{tabular}
\caption{\rev{LVSIS model. Comparison of prediction performance. Spatio-temporal model (1D).}}
\label{fig:test-lvsis-1d}
\end{figure}

\begin{figure}[h!]
\centering
\setlength{\tabcolsep}{4pt}
\renewcommand{\arraystretch}{1.2}
\begin{tabular}{c c c c c}
 & 1/3 & 1/2 & 2/3 \\
\centering\rotatebox{90}{Dynamics} & 
\includegraphics[width=0.3\textwidth]{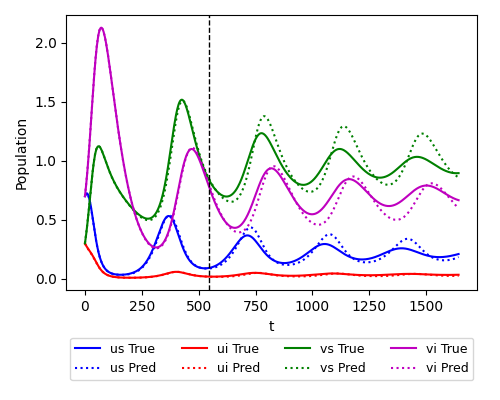} &
\includegraphics[width=0.3\textwidth]{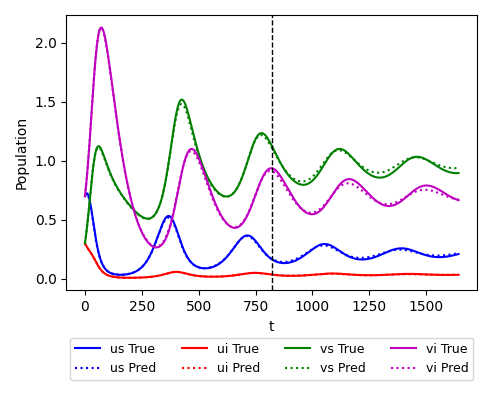} &
\includegraphics[width=0.3\textwidth]{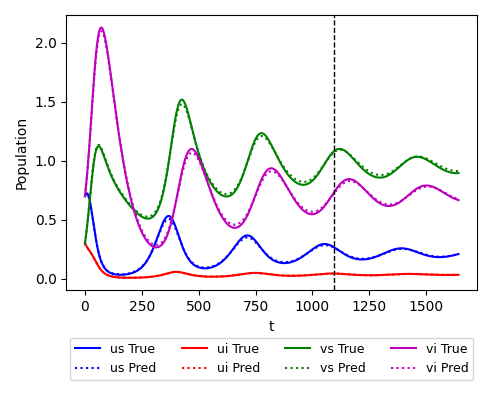} 
\\
\centering\rotatebox{90}{Prediction} & 
\includegraphics[width=0.3\textwidth]{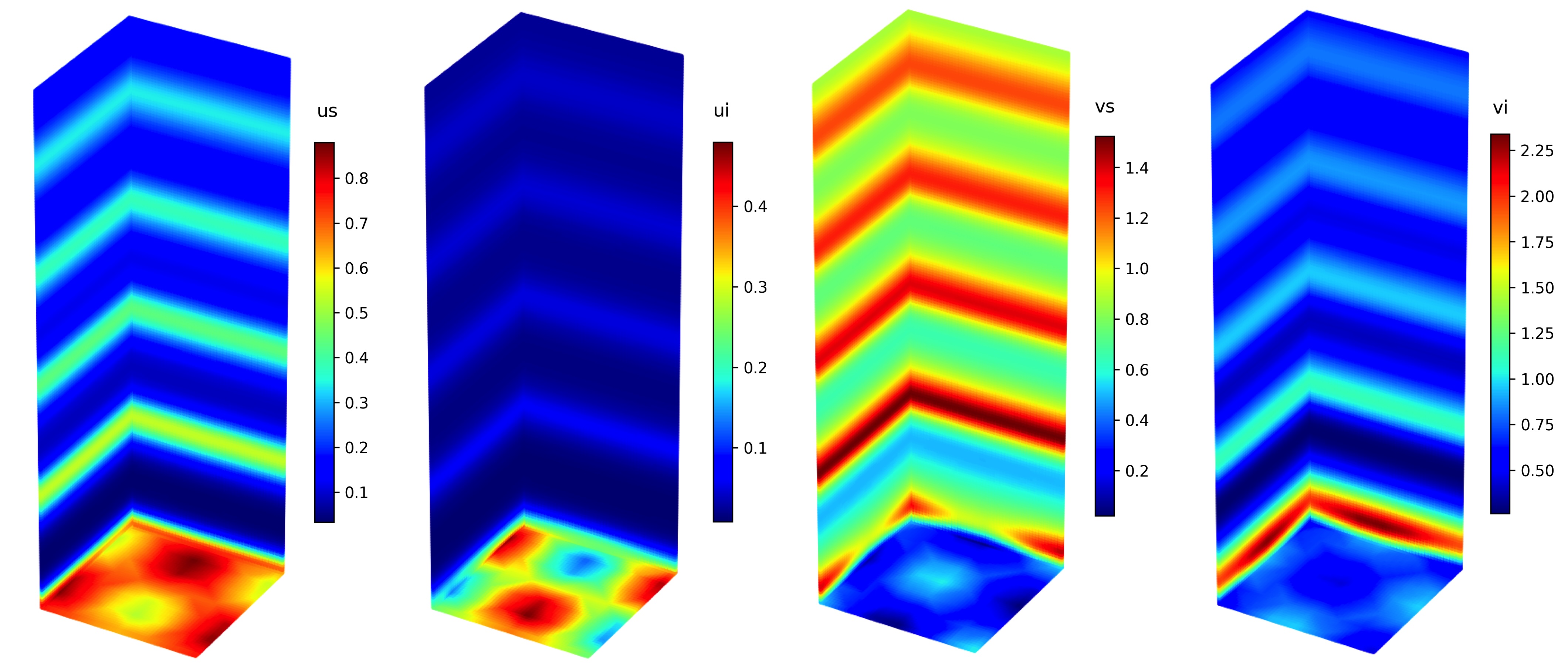} &
\includegraphics[width=0.3\textwidth]{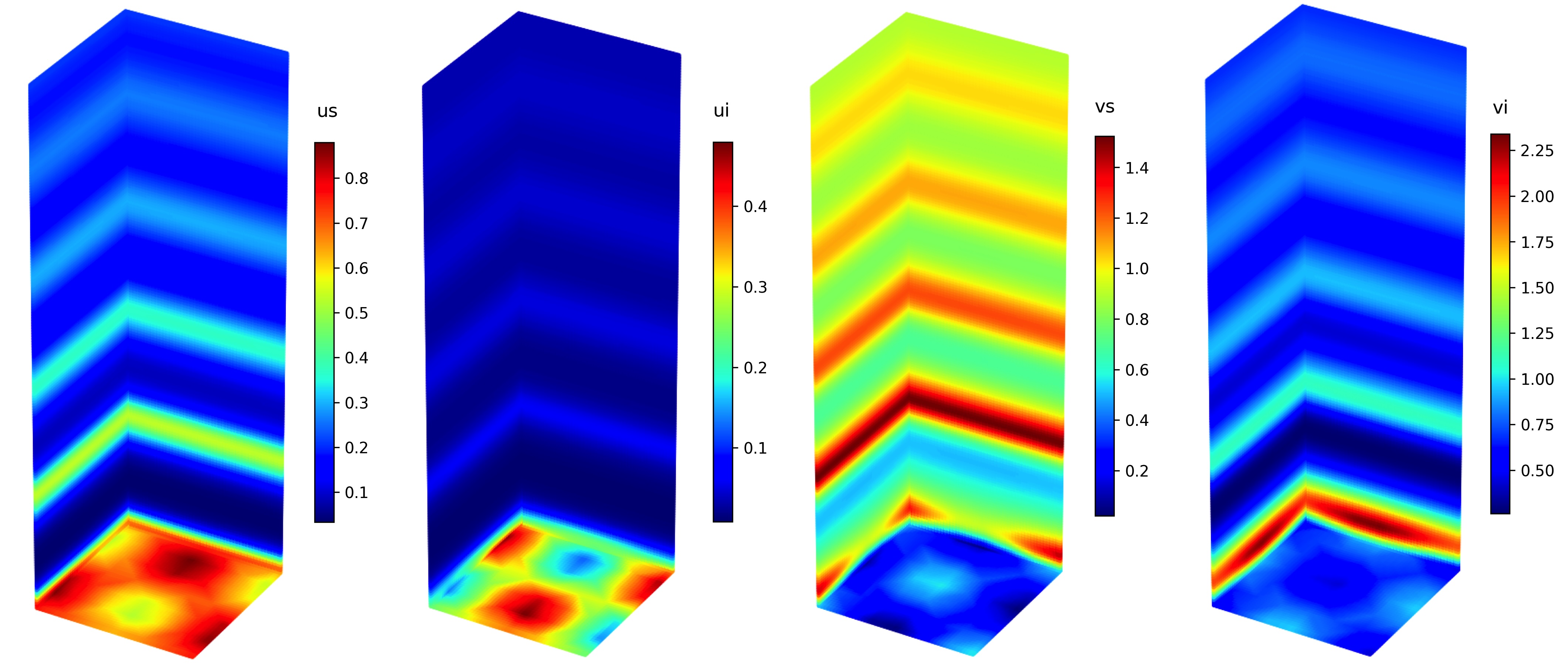} &
\includegraphics[width=0.3\textwidth]{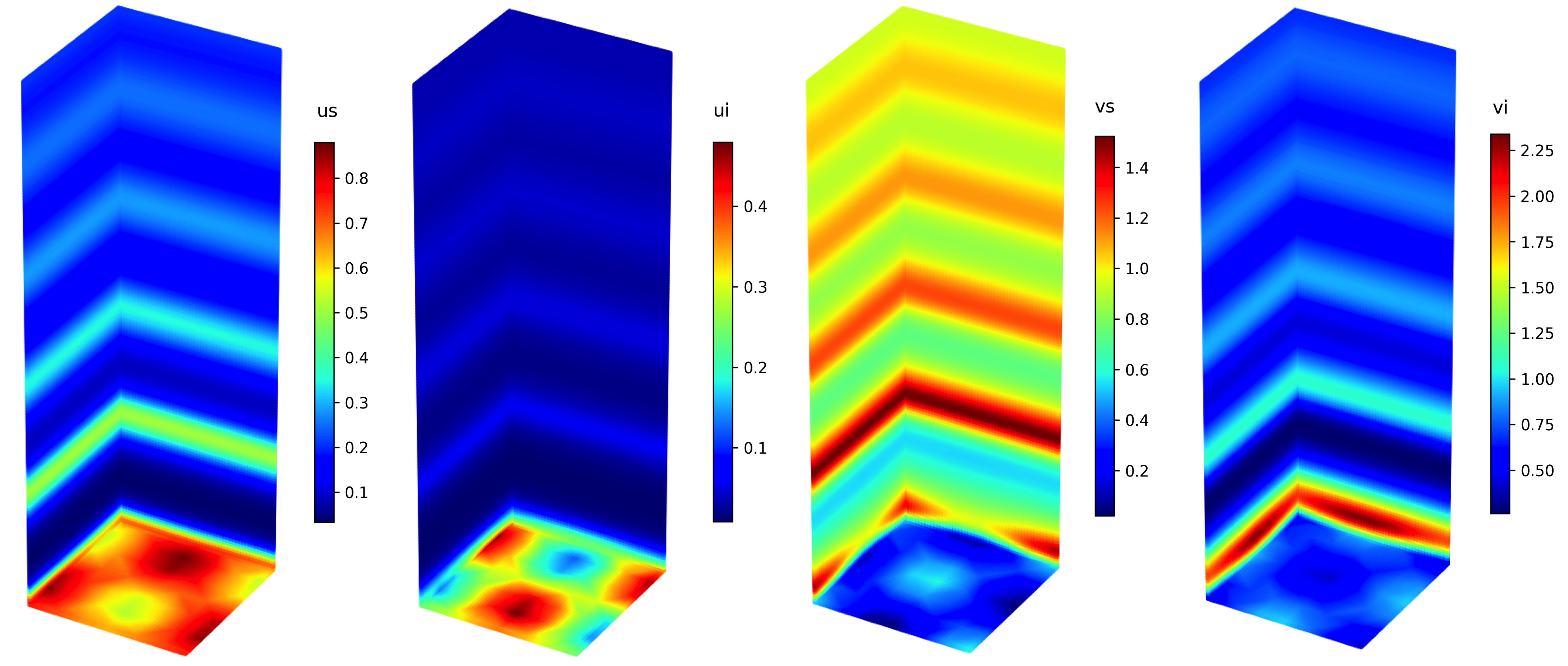} 
\\
\end{tabular}
\caption{\rev{LVSIS model. Comparison of prediction performance. Spatio-temporal model (2D).}}
\label{fig:test-lvsis-2d}
\end{figure}

\begin{table}[h!]
\centering
\small
\begin{tabular}{|l|cc|cc|cc|}
\hline
 & \multicolumn{2}{c|}{1/3} & \multicolumn{2}{c|}{1/2} & \multicolumn{2}{c|}{2/3} \\
Grid & Train & Full & Train & Full & Train & Full \\
\hline
\multicolumn{7}{|c|}{SIR} \\
\hline
1D & 2.472 $\pm$ 0.734 & 54.196 $\pm$ 16.250 & 2.369 $\pm$ 1.193 & 23.062 $\pm$ 15.272 & 1.797 $\pm$ 1.162 & 1.592 $\pm$ 0.662 \\
2D & 3.310 $\pm$ 2.641 & 75.934 $\pm$ 9.114 & 2.509 $\pm$ 1.199 & 23.129 $\pm$ 17.199 & 1.903 $\pm$ 0.551 & 1.577 $\pm$ 0.346 \\
\hline
\multicolumn{7}{|c|}{LV} \\
\hline
1D & 0.258 $\pm$ 0.043 & 0.692 $\pm$ 0.329 & 0.542 $\pm$ 0.071 & 0.501 $\pm$ 0.105 & 0.518 $\pm$ 0.387 & 0.493 $\pm$ 0.404 \\
2D & 0.273 $\pm$ 0.048 & 0.775 $\pm$ 0.085 & 0.417 $\pm$ 0.030 & 0.532 $\pm$ 0.158 & 0.619 $\pm$ 0.203 & 0.569 $\pm$ 0.153 \\
\hline
\multicolumn{7}{|c|}{LVSIS} \\
\hline
1D & 1.372 $\pm$ 0.499 & 14.013 $\pm$ 0.925 & 2.976 $\pm$ 1.584 & 3.485 $\pm$ 1.681 & 2.783 $\pm$ 1.199 & 2.411 $\pm$ 1.023 \\
2D & 2.219 $\pm$ 0.994 & 14.658 $\pm$ 1.651 & 1.647 $\pm$ 0.244 & 1.719 $\pm$ 0.341 & 2.106 $\pm$ 0.713 & 2.046 $\pm$ 0.524 \\
\hline
\end{tabular}
\caption{\rev{Training-interval and full-trajectory relative errors (\%) for the full-space reaction--diffusion NODE in one and two spatial dimensions at training proportions $p=1/3$, $1/2$, and $2/3$. The spatial grids contain $10$ nodes in 1D and $7\times7$ nodes in 2D. Entries report mean $\pm$ standard deviation over random initializations.}}
\label{tab:reaction-diffusion-full-errors}
\end{table}

\begin{table}[h!]
\centering
\small
\begin{tabular}{|l|c|c|c|}
\hline
Grid & 1/3 & 1/2 & 2/3 \\
\hline
\multicolumn{4}{|c|}{SIR} \\
\hline
1D & 77.750 $\pm$ 1.405 & 75.082 $\pm$ 0.691 & 77.944 $\pm$ 2.623 \\
2D & 246.891 $\pm$ 3.622 & 251.145 $\pm$ 4.011 & 247.471 $\pm$ 12.729 \\
\hline
\multicolumn{4}{|c|}{LV} \\
\hline
1D & 76.459 $\pm$ 0.829 & 76.992 $\pm$ 1.716 & 76.259 $\pm$ 1.458 \\
2D & 252.611 $\pm$ 9.505 & 259.221 $\pm$ 9.854 & 247.814 $\pm$ 6.874 \\
\hline
\multicolumn{4}{|c|}{LVSIS} \\
\hline
1D & 75.346 $\pm$ 0.433 & 75.848 $\pm$ 1.353 & 77.566 $\pm$ 1.281 \\
2D & 235.217 $\pm$ 0.744 & 245.627 $\pm$ 9.148 & 248.480 $\pm$ 1.857 \\
\hline
\end{tabular}
\caption{\rev{Wall-clock training time (s) for the full-space reaction--diffusion NODE in one and two spatial dimensions at training proportions $p=1/3$, $1/2$, and $2/3$. The spatial grids contain $10$ nodes in 1D and $7\times7$ nodes in 2D. Entries report mean $\pm$ standard deviation over random initializations.}}
\label{tab:reaction-diffusion-full-time}
\end{table}

\rev{Tables~\ref{tab:reaction-diffusion-full-errors} and~\ref{tab:reaction-diffusion-full-time} summarize the accuracy and cost of learning the local reaction term in the full spatial state. The accuracy depends more strongly on the available training horizon and on the underlying dynamics than on spatial dimension. For SIR, the full-trajectory error decreases from $54.20\%$ in 1D and $75.93\%$ in 2D at $p=1/3$ to approximately $1.6\%$ in both dimensions at $p=2/3$. LVSIS exhibits the same extrapolation challenge at the shortest horizon, with errors near $14\%$, but these decrease to $2.41\%$ in 1D and $2.05\%$ in 2D at $p=2/3$. In contrast, the oscillatory LV dynamics remain below $0.8\%$ for every grid and training proportion. Table~\ref{tab:reaction-diffusion-full-time} shows that the training cost is nearly insensitive to the dynamical system and to $p$, but increases from approximately $75$--$78$~s in 1D to $235$--$259$~s in 2D because each ODE rollout evolves a larger spatial state.}

\begin{figure}[h!]
\centering
\includegraphics[width=0.48\linewidth]{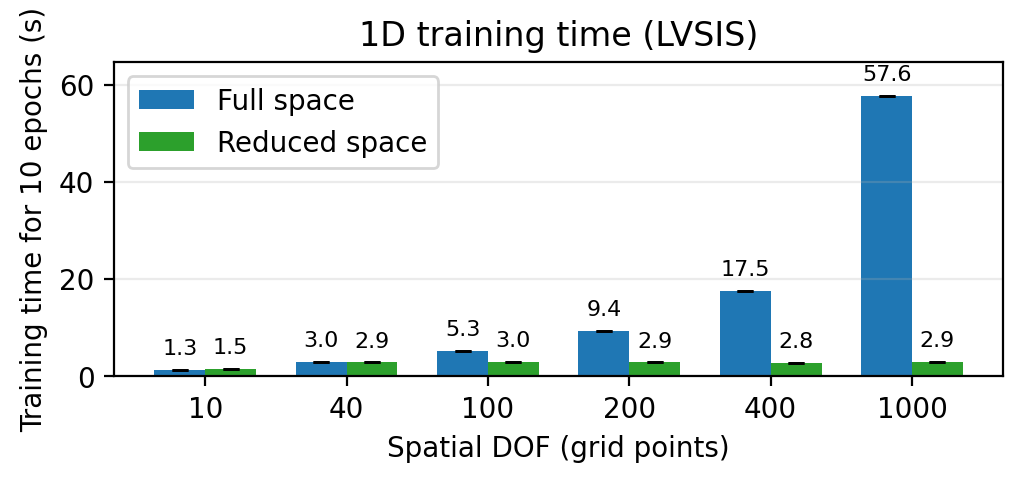}
\includegraphics[width=0.48\linewidth]{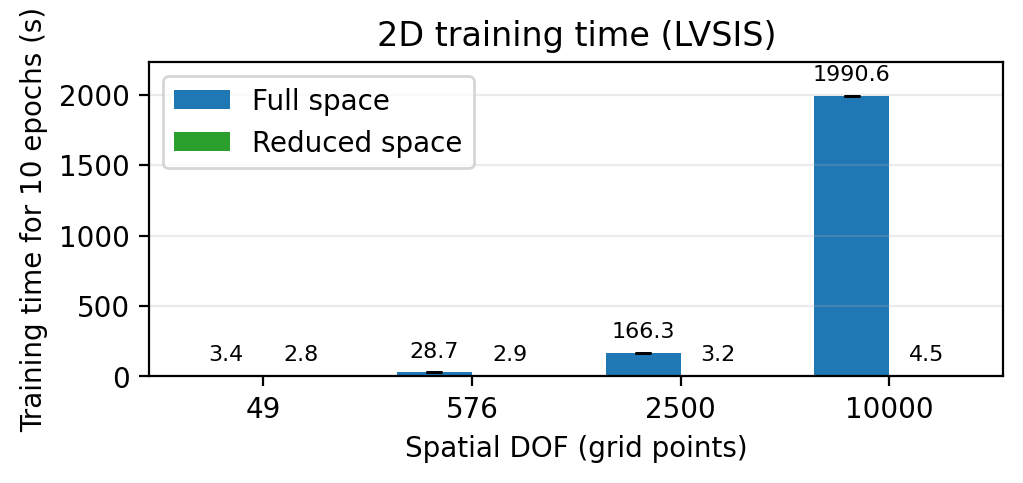}
\caption{\rev{Training-time comparison over $10$ epochs between full-space and reduced space reaction--diffusion learning in one and two spatial dimensions. The reduced space model uses configuration $(r,q)=(4,30)$, where $r$ is the state-basis rank and $q$ is the number of interpolation nodes.}}
\label{fig:reaction-diffusion-training-time}
\end{figure}

\begin{table}[h!]
\centering
\small
\begin{tabular}{|l|cc|cc|cc|}
\hline
 & \multicolumn{2}{c|}{1/3} & \multicolumn{2}{c|}{1/2} & \multicolumn{2}{c|}{2/3} \\
Grid & Train & Full & Train & Full & Train & Full \\
\hline
\multicolumn{7}{|c|}{SIR} \\
\hline
1D & 3.110 $\pm$ 1.119 & 58.019 $\pm$ 15.089 & 2.272 $\pm$ 1.321 & 22.034 $\pm$ 15.259 & 3.019 $\pm$ 1.493 & 2.458 $\pm$ 0.972 \\
2D & 2.335 $\pm$ 1.209 & 68.259 $\pm$ 8.795 & 2.648 $\pm$ 0.656 & 22.930 $\pm$ 16.875 & 2.156 $\pm$ 1.071 & 1.733 $\pm$ 0.774 \\
\hline
\multicolumn{7}{|c|}{LV} \\
\hline
1D & 0.270 $\pm$ 0.025 & 0.664 $\pm$ 0.271 & 0.383 $\pm$ 0.083 & 0.448 $\pm$ 0.187 & 0.576 $\pm$ 0.229 & 0.534 $\pm$ 0.248 \\
2D & 0.550 $\pm$ 0.175 & 0.923 $\pm$ 0.372 & 0.411 $\pm$ 0.047 & 0.501 $\pm$ 0.224 & 0.621 $\pm$ 0.194 & 0.570 $\pm$ 0.157 \\
\hline
\multicolumn{7}{|c|}{LVSIS} \\
\hline
1D & 1.234 $\pm$ 0.198 & 14.012 $\pm$ 0.889 & 2.933 $\pm$ 0.987 & 2.836 $\pm$ 0.478 & 3.076 $\pm$ 1.559 & 2.851 $\pm$ 1.505 \\
2D & 1.535 $\pm$ 0.919 & 14.124 $\pm$ 1.340 & 2.786 $\pm$ 0.400 & 2.754 $\pm$ 0.216 & 1.917 $\pm$ 0.446 & 1.908 $\pm$ 0.378 \\
\hline
\end{tabular}
\caption{\rev{Training-interval and full-trajectory relative errors (\%) for the reduced space reaction--diffusion NODE with configuration $(r,q)=(7,30)$ at training proportions $p=1/3$, $1/2$, and $2/3$. The spatial grids contain $101$ nodes in 1D and $51\times51$ nodes in 2D. Entries report mean $\pm$ standard deviation over random initializations.}}
\label{tab:reaction-diffusion-pod-errors}
\end{table}

\begin{table}[h!]
\centering
\small
\begin{tabular}{|l|c|c|c|}
\hline
Grid & 1/3 & 1/2 & 2/3 \\
\hline
\multicolumn{4}{|c|}{SIR} \\
\hline
1D & 183.042 $\pm$ 2.869 & 180.319 $\pm$ 3.629 & 188.976 $\pm$ 2.906 \\
2D & 240.495 $\pm$ 5.021 & 235.752 $\pm$ 7.236 & 231.131 $\pm$ 3.323 \\
\hline
\multicolumn{4}{|c|}{LV} \\
\hline
1D & 183.961 $\pm$ 1.006 & 189.799 $\pm$ 2.475 & 183.241 $\pm$ 4.253 \\
2D & 211.356 $\pm$ 2.681 & 205.476 $\pm$ 1.985 & 207.316 $\pm$ 1.757 \\
\hline
\multicolumn{4}{|c|}{LVSIS} \\
\hline
1D & 190.371 $\pm$ 8.217 & 188.885 $\pm$ 3.061 & 197.294 $\pm$ 18.602 \\
2D & 214.732 $\pm$ 4.443 & 212.004 $\pm$ 2.267 & 216.292 $\pm$ 0.468 \\
\hline
\end{tabular}
\caption{\rev{Wall-clock training time (s) for the reduced space reaction--diffusion NODE with configuration $(r,q)=(7,30)$ at training proportions $p=1/3$, $1/2$, and $2/3$. The spatial grids contain $101$ nodes in 1D and $51\times51$ nodes in 2D. Entries report mean $\pm$ standard deviation over random initializations.}}
\label{tab:reaction-diffusion-pod-time}
\end{table}

\rev{Figure~\ref{fig:reaction-diffusion-training-time} isolates the spatial-scaling effect for LVSIS using the POD--QDEIM configuration $(r,q)=(4,30)$. On the smallest grids, reduction provides small benefit since basis projection and QDEIM sampling introduce fixed overhead. As the number of spatial degrees of freedom increases, full-space cost grows rapidly whereas the reduced cost remains nearly constant. In 1D, the $10$-epoch time increases from $1.3$~s at $10$ grid points to $57.6$~s at $1000$ grid points, while POD--QDEIM remains near $2.8$--$3.0$~s, yielding an approximately $20$-fold speedup at the largest grid. In 2D, full-space time rises from $3.4$~s at $49$ degrees of freedom to $1990.6$~s at $10{,}000$, while the reduced time increases only from $2.8$ to $4.5$~s. }

\rev{Tables~\ref{tab:reaction-diffusion-pod-errors} and~\ref{tab:reaction-diffusion-pod-time} assess the POD--QDEIM configuration $(r,q)=(7,30)$ on the finer $101$-node 1D and $51\times51$-node 2D grids. LV full-trajectory errors remain below $0.93\%$. For SIR and LVSIS, the reduced model retains the full-space dependence on training horizon: errors are large at $p=1/3$ but fall to $1.73$--$2.46\%$ for SIR and $1.91$--$2.85\%$ for LVSIS at $p=2/3$. The corresponding training times remain within approximately $180$--$197$~s in 1D and $205$--$240$~s in 2D across the three systems, indicating substantially weaker spatial-size sensitivity than a full-space solve on comparably fine grids.}

\begin{table}[h!]
\centering
\small
\begin{tabular}{|l|cc|cc|cc|}
\hline
 & \multicolumn{2}{c|}{Full space} & \multicolumn{2}{c|}{\shortstack{Reduced space \\$(r=4,q=30)$}} & \multicolumn{2}{c|}{\shortstack{Reduced space \\$(r=7,q=30)$}} \\
Grid & Train & Full & Train & Full & Train & Full \\
\hline
\multicolumn{7}{|c|}{LVSIS, noise-free} \\
\hline
1D & 2.125 $\pm$ 0.414 & 2.150 $\pm$ 0.123 & 1.796 $\pm$ 0.561 & 2.096 $\pm$ 0.721 & 3.158 $\pm$ 0.381 & 2.933 $\pm$ 0.243 \\
2D & 2.261 $\pm$ 0.907 & 2.939 $\pm$ 0.821 & 3.066 $\pm$ 1.542 & 3.000 $\pm$ 1.074 & 2.222 $\pm$ 0.198 & 2.061 $\pm$ 0.172 \\
\hline
\multicolumn{7}{|c|}{LVSIS, noise-free + SINDy} \\
\hline
1D & 3.740 $\pm$ 1.428 & 3.118 $\pm$ 1.119 & 3.831 $\pm$ 1.632 & 3.168 $\pm$ 1.403 &- & -\\
2D & 3.350 $\pm$ 0.280 & 3.160 $\pm$ 0.502 & 5.543 $\pm$ 3.846 & 4.508 $\pm$ 2.772 &- & - \\
\hline
\multicolumn{7}{|c|}{LVSIS, noisy ($\sigma=0.01$)} \\
\hline
1D & 6.246 $\pm$ 4.095 & 6.751 $\pm$ 3.994 & 3.901 $\pm$ 0.586 & 4.982 $\pm$ 1.220 & 11.581 $\pm$ 7.101 & 101.078 $\pm$ 131.494 \\
2D & 4.185 $\pm$ 1.837 & 4.098 $\pm$ 1.110 & 3.702 $\pm$ 1.163 & 3.578 $\pm$ 0.558 & 20.717 $\pm$ 22.670 & 76.135 $\pm$ 87.970 \\
\hline
\multicolumn{7}{|c|}{LVSIS, noisy ($\sigma=0.01$) + SINDy} \\
\hline
1D & 9.640 $\pm$ 5.477 & 9.399 $\pm$ 4.830 & 5.722 $\pm$ 2.104 & 5.323 $\pm$ 1.750 &- & - \\
2D & 5.022 $\pm$ 0.872 & 4.530 $\pm$ 0.734 & 4.548 $\pm$ 0.000 & 3.917 $\pm$ 0.000 &- & - \\
\hline
\end{tabular}
\caption{\rev{Training-interval and full-trajectory relative errors (\%) for the full-space LVSIS model and reduced space configurations $(r,q)=(4,30)$ and $(7,30)$ at training proportion $p=1/2$. Results are shown for noise-free data and additive Gaussian noise at $\sigma=0.01$, before and after the optional SINDy identification stage. Entries report mean $\pm$ standard deviation over random initializations.}}
\label{tab:reaction-diffusion-pod-ranks}
\end{table}

\rev{Tables~\ref{tab:reaction-diffusion-pod-ranks} and~\ref{tab:reaction-diffusion-pod-ranks-time} compare the full-space model with POD--QDEIM configurations $(r,q)=(4,30)$ and $(7,30)$ for LVSIS at $p=1/2$. For noise-free data, reduction does not materially degrade accuracy: in 1D, full-trajectory errors are $2.15\%$, $2.10\%$, and $2.93\%$ for the full-space, $(4,30)$, and $(7,30)$ configurations, respectively; in 2D, $(7,30)$ gives the smallest error, $2.06\%$, compared with $2.94\%$ for the full-space model. For additive Gaussian noise with $\sigma=0.01$, $(4,30)$ is the most reliable reduced configuration, lowering the full-trajectory error from $6.75\%$ to $4.98\%$ in 1D and from $4.10\%$ to $3.58\%$ in 2D. At the same noise level, $(7,30)$ is unstable, with full errors above $76\%$ and very large standard deviations; its larger errors and poor denoising make the resulting trajectory unsuitable for the SINDy stage. Increasing the reduced dimension therefore does not automatically improve robustness. For a reliable learned trajectory, the optional SINDy stage can then be applied to obtain a sparse, interpretable representation and recover an explicit governing equation.}

\begin{table}[h!]
\centering
\small
\begin{tabular}{|l|c|c|c|}
\hline 
Grid & Full space 
& \shortstack{Reduced space, \\ $(r=4,q=30)$} & 
\shortstack{Reduced space, \\ $(r=7,q=30)$} \\
\hline
\multicolumn{4}{|c|}{LVSIS} \\
\hline
1D & 210.211 $\pm$ 1.835 & 197.622 $\pm$ 3.644 & 193.660 $\pm$ 0.060 \\
2D & 2269.656 $\pm$ 7.408 & 201.607 $\pm$ 2.831 & 200.051 $\pm$ 2.064 \\
\hline
\end{tabular}
\caption{\rev{Wall-clock training time (s) for the noise-free full-space LVSIS model and reduced space configurations $(r,q)=(4,30)$ and $(7,30)$ at training proportion $p=1/2$. Entries report mean $\pm$ standard deviation over random initializations.}}
\label{tab:reaction-diffusion-pod-ranks-time}
\end{table}

\rev{Table~\ref{tab:reaction-diffusion-pod-ranks-time} quantifies the computational benefit on the fine grids of $101$ nodes in 1D and $51\times51$ nodes in 2D. In 1D, training time decreases from $210.21$~s for the full-space model to $197.62$~s for $(4,30)$ and $193.66$~s for $(7,30)$. In 2D, it decreases from $2269.66$~s to approximately $200$~s for either reduced configuration, an eleven-fold speedup while preserving approximately $2$--$3\%$ full-trajectory error on noise-free data. Across Tables~\ref{tab:reaction-diffusion-full-errors}--\ref{tab:reaction-diffusion-pod-ranks-time} and Figure~\ref{fig:reaction-diffusion-training-time}, the benefit of POD--QDEIM is clearest for large spatial states, while the noise results show that $r$ must be selected together with the data quality.}

\section{Conclusion}

\rev{We considered SINDy, Neural ODEs, and KANODEs for learning eco--epidemiological dynamics from time-series data. All three methods captured the coupled nonlinear dynamics, but with different strengths. SINDy recovered sparse governing equations and provided direct interpretability, whereas Neural ODEs and KANODEs offered greater flexibility for complex and noisy dynamics. KANODEs achieved comparable accuracy with fewer trainable parameters, while deeper and wider Neural ODEs were particularly effective for the more challenging LVSIS system.}

\rev{We then extended neural trajectory learning to reaction--diffusion systems in which the diffusion operator is prescribed and only the hidden local reaction is learned. The formulation provides a two-level reduction: POD compresses the spatial state and known diffusion operator, while QDEIM hyper-reduces the NODE reaction evaluation to selected interpolation nodes. Our formulation retains reconstruction and error evaluation on the original grid and provides its largest computational savings on fine two-dimensional grids. The rank study further shows that the reduced dimension must be selected according to data quality: a smaller rank can improve denoising, whereas a larger rank may amplify noisy directions and prevent reliable SINDy equation recovery.}

\rev{We note that the present evaluation uses physics-based synthetic data from canonical models, providing known ground truth and a structured setting in which to isolate the effects of noise, training horizon, spatial resolution, and reduced-space approximation. The spatial framework also assumes that the diffusion operator is known, which may not hold in many applications. For more complex dynamics, the SINDy candidate library may also need to be enlarged beyond quadratic polynomials. These controlled experiments do not replace validation on empirical ecological or epidemic observations. Validation on suitable real-world ecological and epidemic datasets therefore remains an important direction for future work.}

\bibliographystyle{unsrt}
\bibliography{lit}

\clearpage
\appendix
\section{\rev{Supplementary Temporal Results}}
\label{app:temporal-supplement}

\rev{This appendix presents the supplementary temporal results in one continuous sequence, following the organization used previously in the main text. It includes the noise-free trajectory comparisons for the SIR, LV, and LVSIS benchmarks, the corresponding sparse equations identified by SINDy, and the noisy-data trajectory comparisons.}

\begin{figure}[h!]
\centering
\setlength{\tabcolsep}{4pt}
\renewcommand{\arraystretch}{1.2}
\begin{tabular}{c c c c c}
 & {1/3} & {1/2} & {2/3} \\
\centering\rotatebox{90}{SINDy} &
\includegraphics[width=0.3\textwidth]{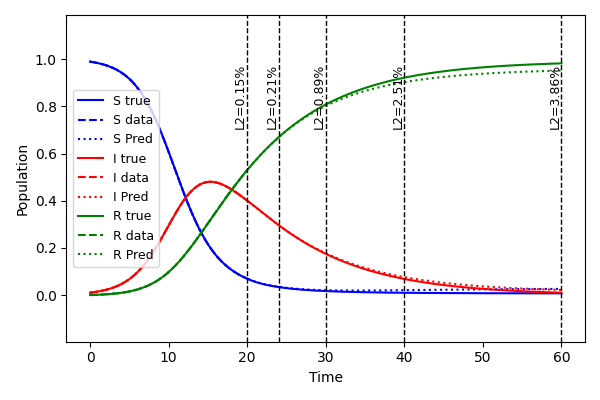} & 
\includegraphics[width=0.3\textwidth]{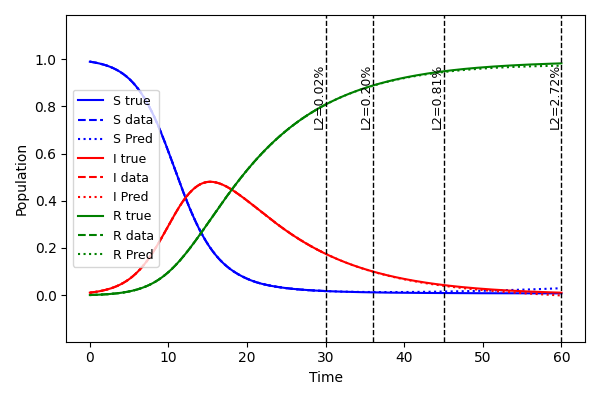} & 
\includegraphics[width=0.3\textwidth]{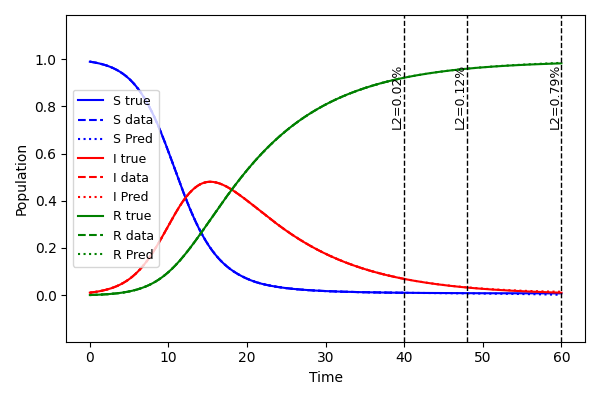} 
\\
\centering\rotatebox{90}{NODEs} & 
\includegraphics[width=0.3\textwidth]{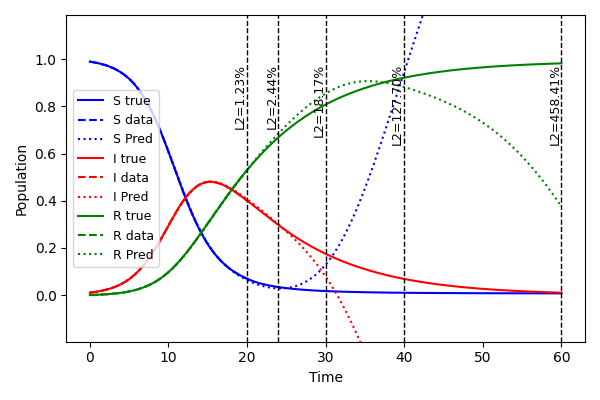} &
\includegraphics[width=0.3\textwidth]{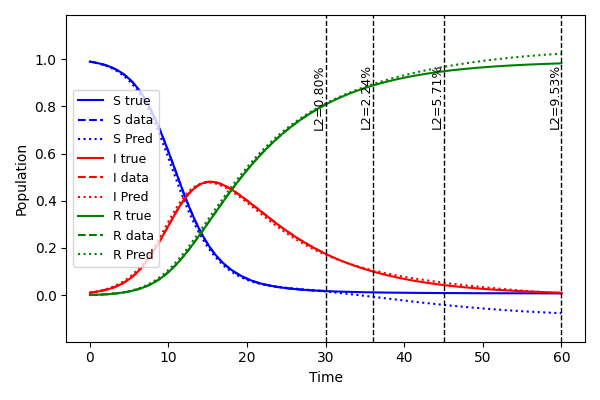} &
\includegraphics[width=0.3\textwidth]{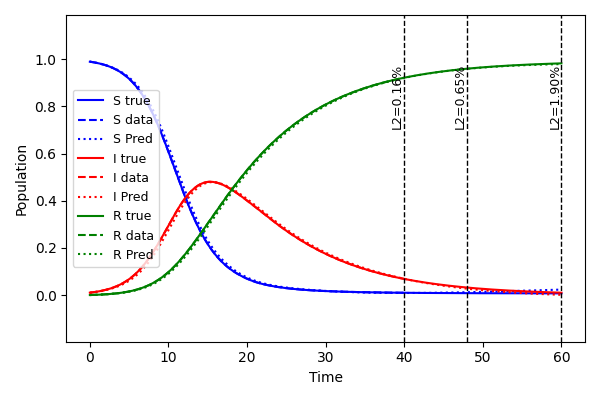} 
\\
\centering\rotatebox{90}{KANODEs} &
\includegraphics[width=0.3\textwidth]{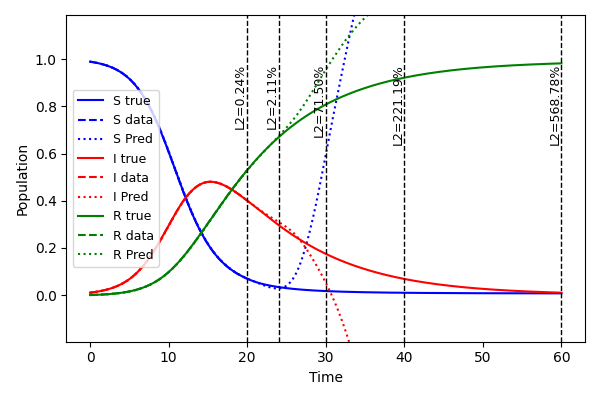} &
\includegraphics[width=0.3\textwidth]{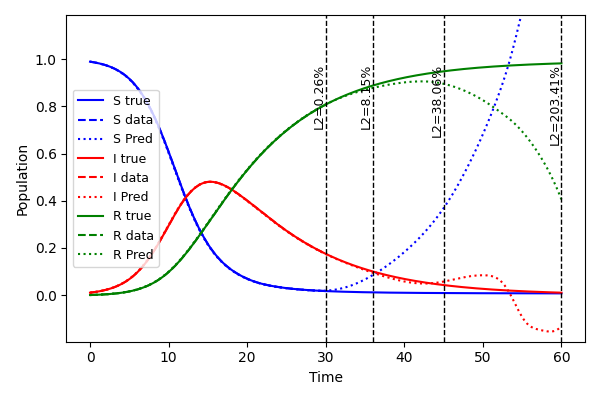} &
\includegraphics[width=0.3\textwidth]{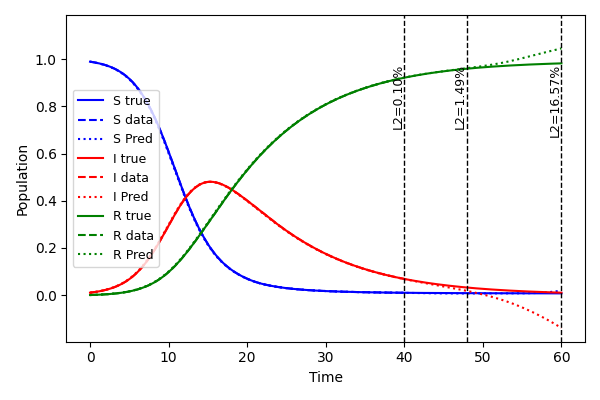} 
\\
\end{tabular}
\caption{\rev{SIR model. Comparison of prediction performance.}}
\label{fig:test-sir1}
\end{figure}

\begin{table}[h!]
\centering
\scriptsize
\begin{tabular}{|c|c|}
\hline
\textbf{} & \textbf{Learned Equations} \\
\hline
\textbf{1/3}  & 
$
\begin{aligned}
u_s' &= -0.065  -0.104 u_s  -0.001 u_i + 0.040 u_r + 0.169 u_s^2  {-0.198 u_s u_i}  -0.075 u_s u_r + 0.116 u_i^2 + 0.081 u_i u_r + 0.034 u_r^2\\
u_i' &= 0.024  + 0.087 u_s \hr{+ 0.016 u_i}  -0.079 u_r  -0.112 u_s^2 \hr{+ 0.245 u_s u_i}  -0.047 u_s u_r \hb{-0.141 u_i^2} - 0.088 u_i u_r + 0.055 u_r^2\\
u_r' &= 0.041  + 0.017 u_s \hr{- 0.016 u_i} + 0.039 u_r - 0.057 u_s^2 - 0.048 u_s u_i + 0.122 u_s u_r \hb{+ 0.025 u_i^2} + 0.007 u_i u_r - 0.089 u_r^2
\end{aligned}
$
\\  \hline
\textbf{1/2}  & 
$
\begin{aligned}
u_s' &= -0.046  - 0.068 u_s - 0.078 u_i + 0.100 u_r + 0.115 u_s^2 \hr{- 0.268 u_s u_i} + 0.085 u_s u_r + 0.123 u_i^2 + 0.067 u_i u_r - 0.052 u_r^2\\
u_i' &= 0.019  + 0.076 u_s \hr{+ 0.039 u_i} - 0.096 u_r - 0.096 u_s^2 \hr{+ 0.266 u_s u_i} - 0.094 u_s u_r \hb{- 0.148 u_i^2} - 0.079 u_i u_r + 0.077 u_r^2\\
u_r' &= 0.027  - 0.008 u_s \hr{+ 0.039 u_i} - 0.004 u_r - 0.019 u_s^2 + 0.002 u_s u_i + 0.009 u_s u_r \hb{+ 0.025 u_i^2} + 0.012 u_i u_r - 0.024 u_r^2
\end{aligned}
$
\\  \hline
\textbf{2/3}  & 
$
\begin{aligned}
u_s' &= -0.040  - 0.056 u_s - 0.105 u_i + 0.121 u_r + 0.096 u_s^2 \hr{- 0.292 u_s u_i} + 0.141 u_s u_r + 0.125 u_i^2 + 0.062 u_i u_r - 0.082 u_r^2\\
u_i' &= 0.016  + 0.070 u_s \hr{+ 0.053 u_i} - 0.107 u_r - 0.086 u_s^2 \hr{+ 0.278 u_s u_i} - 0.122 u_s u_r \hb{- 0.149 u_i^2} - 0.076 u_i u_r + 0.092 u_r^2\\
u_r' &= 0.024  - 0.014 u_s \hr{+ 0.052 u_i} - 0.014 u_r - 0.010 u_s^2 + 0.014 u_s u_i - 0.019 u_s u_r \hb{+ 0.023 u_i^2} + 0.014 u_i u_r - 0.010 u_r^2
\end{aligned}
$
\\ \hline
\textbf{} & \textbf{Original Equations} \\
\hline
  & 
$
\begin{aligned}
u_s' &=  \hr{-0.5 u_s u_i} \\
u_i' &= \hr{- 0.1 u_i + 0.5 u_s u_i}\\
u_r' &= \hr{0.1 u_i}
\end{aligned}
$
\\  \hline
\end{tabular}
\caption{\rev{SINDy equations identified for the SIR temporal system at training proportions $p=1/3$, $1/2$, and $2/3$; the original equations are included for comparison.}}
\label{tab:sir}
\end{table}

\begin{figure}[h!]
\centering
\setlength{\tabcolsep}{4pt}
\renewcommand{\arraystretch}{1.2}
\begin{tabular}{c c c c c}
 & {1/3} & {1/2} & {2/3} \\
\centering\rotatebox{90}{SINDy} &
\includegraphics[width=0.3\textwidth]{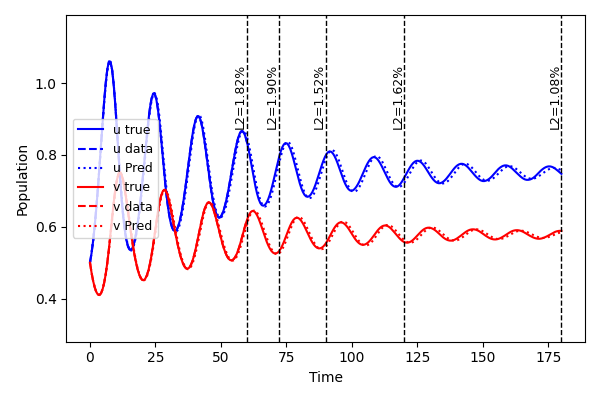} & 
\includegraphics[width=0.3\textwidth]{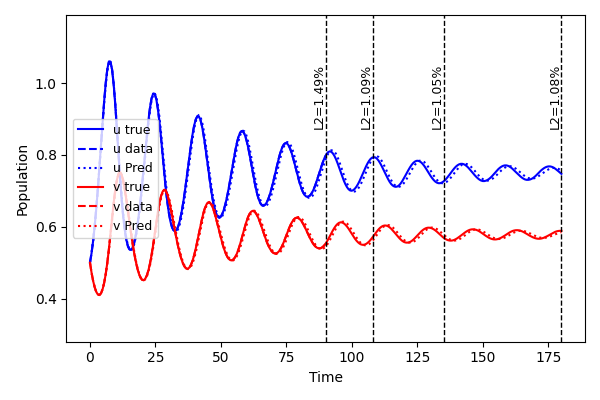} & 
\includegraphics[width=0.3\textwidth]{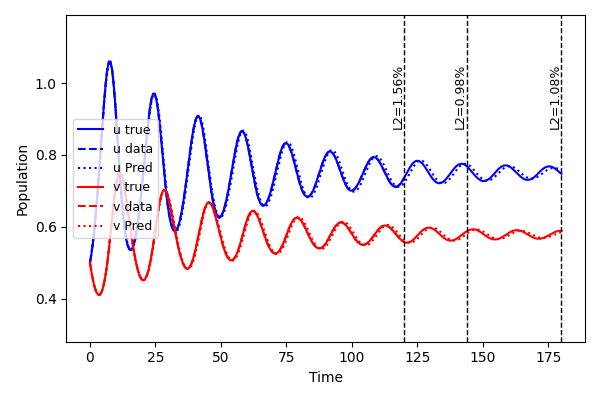} 
\\
\centering\rotatebox{90}{NODEs} & 
\includegraphics[width=0.3\textwidth]{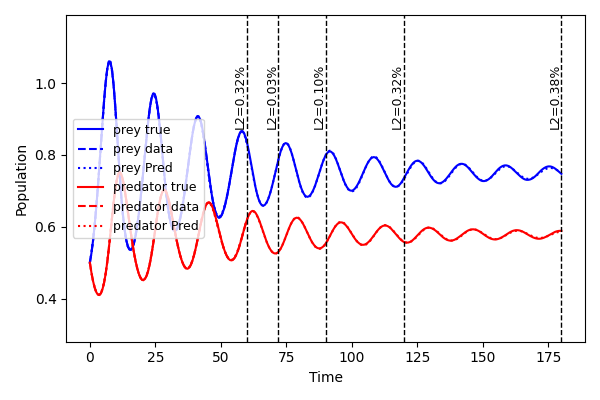} &
\includegraphics[width=0.3\textwidth]{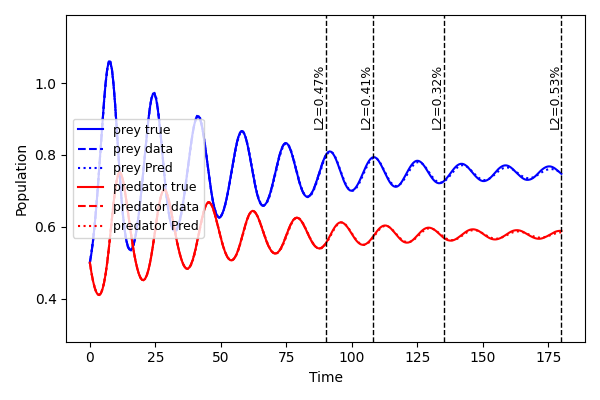} &
\includegraphics[width=0.3\textwidth]{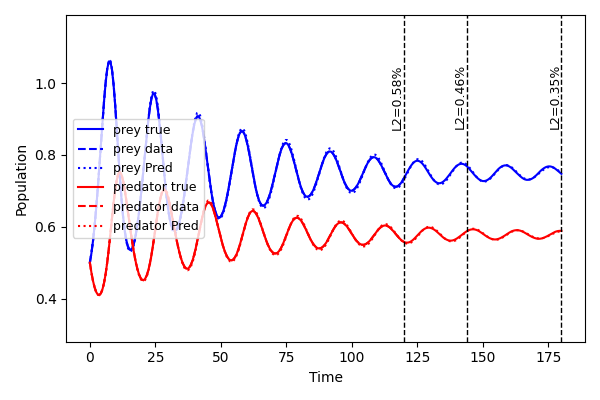} 
\\
\centering\rotatebox{90}{KANODEs} &
\includegraphics[width=0.3\textwidth]{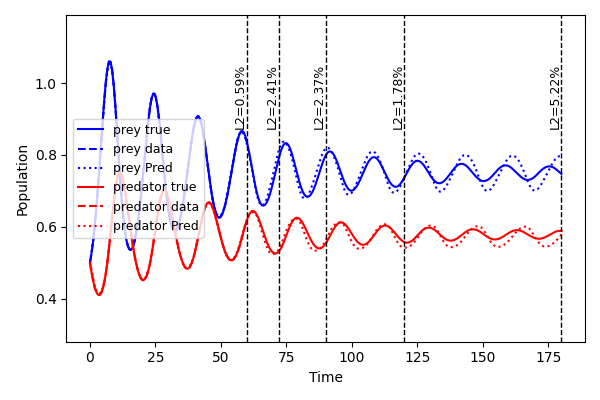} &
\includegraphics[width=0.3\textwidth]{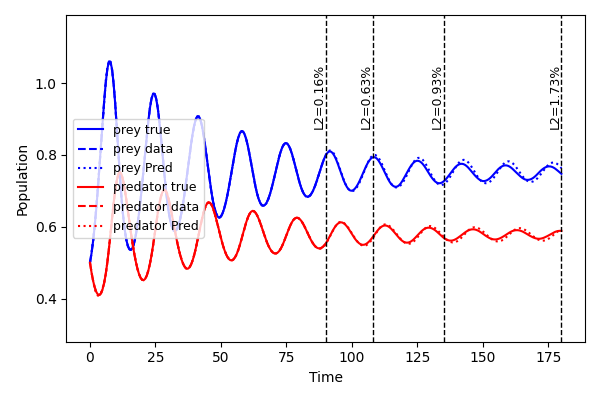} &
\includegraphics[width=0.3\textwidth]{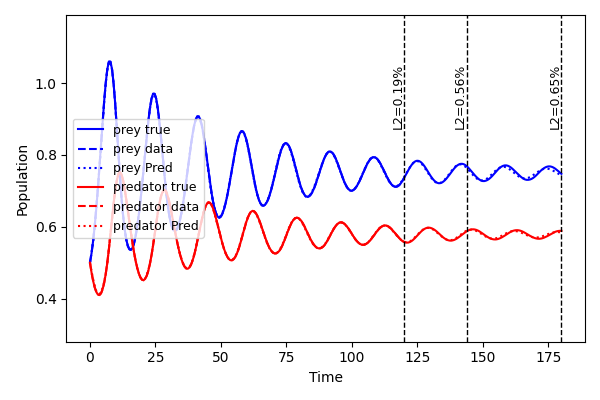} 
\\
\end{tabular}
\caption{\rev{LV model. Comparison of prediction performance.}}
\label{fig:test-lv1}
\end{figure}

\begin{table}[h!]
\centering
\scriptsize
\begin{tabular}{|c|c|}
\hline
\textbf{} & \textbf{Learned Equations} \\
\hline
\textbf{1/3}  & 
$
\begin{aligned}
u' &= 0.002  \hr{+ 0.472 u} + 0.028 v \hr{- 0.043 u^2 - 0.769 u v} - 0.043 v^2\\
v' &= -0.001  + 0.017 u \hr{- 0.320 v} - 0.007 u^2 \hr{+ 0.386 u v} + 0.028 v^2
\end{aligned}
$
\\  \hline
\textbf{1/2}  & 
$
\begin{aligned}
u' &= 0.001  \hr{+ 0.472 u} + 0.031 v \hr{- 0.043 u^2 - 0.769 u v} - 0.046 v^2\\
v' &= 0.017 u \hr{- 0.322 v} - 0.006 u^2 + \hr{0.386 u v} + 0.030 v^2
\end{aligned}
$
\\  \hline
\textbf{2/3}  & 
$
\begin{aligned}
u' &= 0.002  \hr{+ 0.471 u} + 0.029 v \hr{- 0.042 u^2 - 0.769 u v} - 0.045 v^2\\
v' &= 0.017 u \hr{- 0.322 v} - 0.006 u^2 \hr{+ 0.386 u v} + 0.030 v^2
\end{aligned}
$
\\ \hline
\textbf{} & \textbf{Original Equations} \\
\hline
  & 
$
\begin{aligned}
u' &=  \hr{0.5 u - 0.05 u^2 - 0.8 uv} \\
v' &= \hr{-0.3v + 0.4 uv}
\end{aligned}
$
\\  \hline
\end{tabular}
\caption{\rev{SINDy equations identified for the LV temporal system at training proportions $p=1/3$, $1/2$, and $2/3$; the original equations are included for comparison.}}
\label{tab:lv}
\end{table}

\begin{figure}[h!]
\centering
\setlength{\tabcolsep}{4pt}
\renewcommand{\arraystretch}{1.2}
\begin{tabular}{c c c c c}
 & {1/3} & {1/2} & {2/3} \\
\centering\rotatebox{90}{SINDy} &
\includegraphics[width=0.3\textwidth]{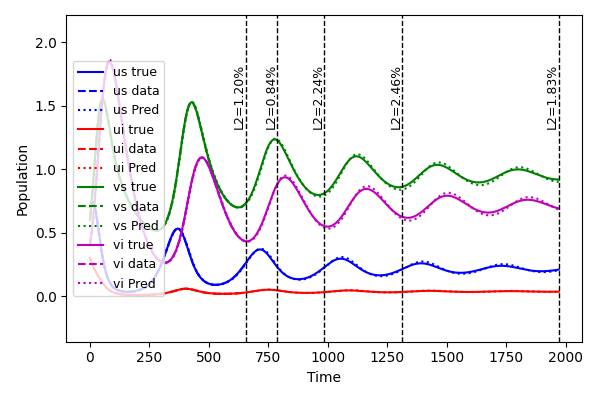} & 
\includegraphics[width=0.3\textwidth]{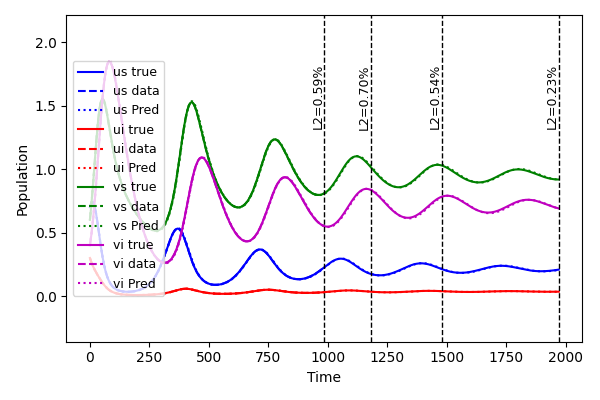} & 
\includegraphics[width=0.3\textwidth]{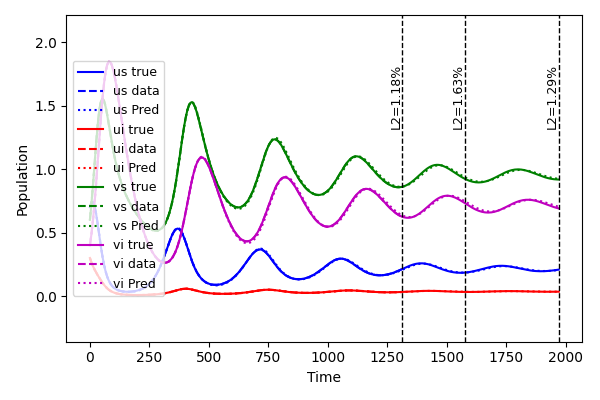} 
\\
\centering\rotatebox{90}{NODEs} & 
\includegraphics[width=0.3\textwidth]{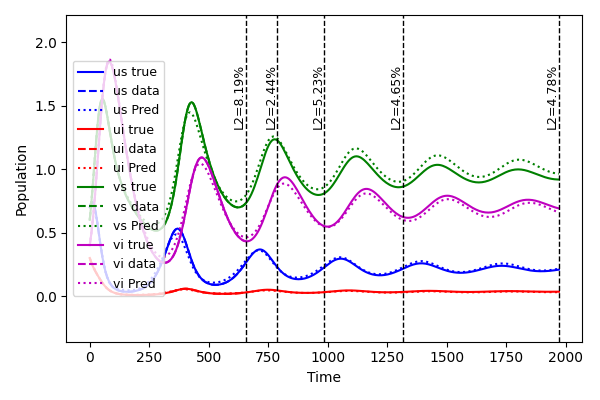} &
\includegraphics[width=0.3\textwidth]{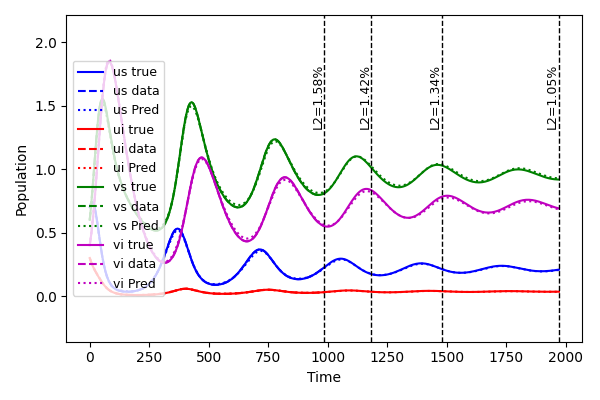} &
\includegraphics[width=0.3\textwidth]{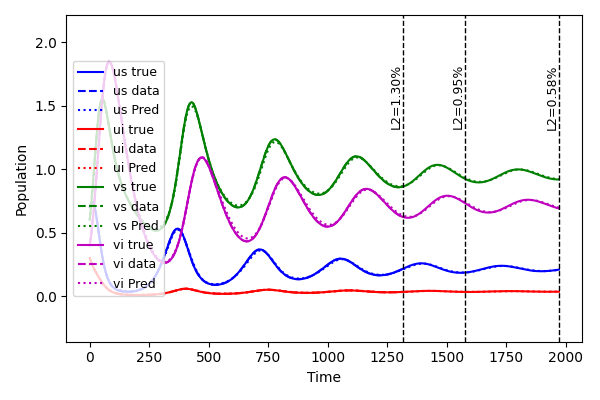} 
\\
\centering\rotatebox{90}{KANODEs} &
\includegraphics[width=0.3\textwidth]{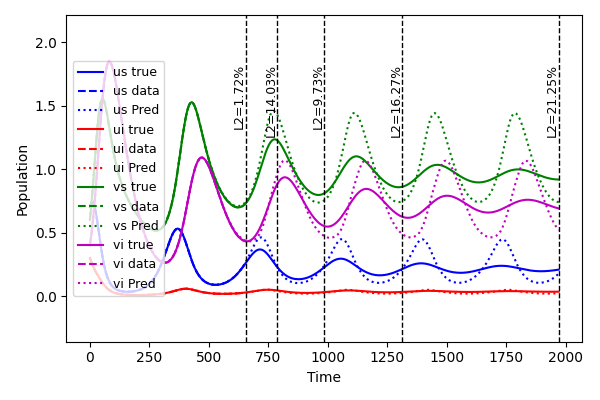} &
\includegraphics[width=0.3\textwidth]{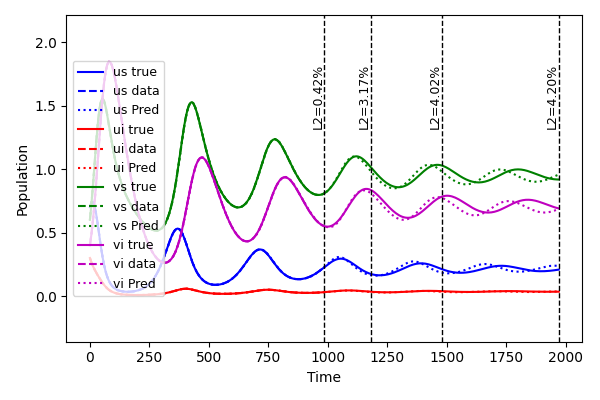} &
\includegraphics[width=0.3\textwidth]{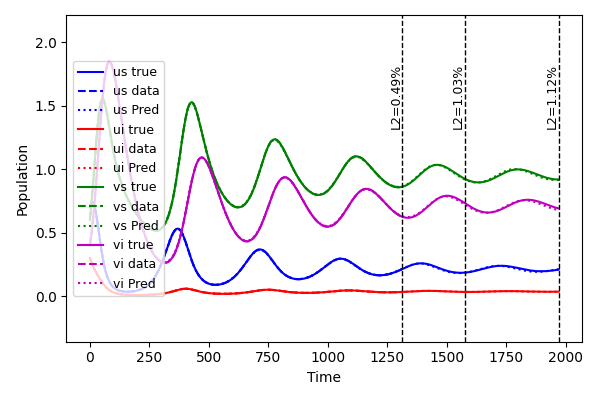} 
\\
\end{tabular}
\caption{\rev{LVSIS model. Comparison of prediction performance.}}
\label{fig:test-lvsis1}
\end{figure}

\begin{table}[h!]
\centering
\scriptsize
\begin{tabular}{|c|c|}
\hline
\textbf{} & \textbf{Learned Equations} \\
\hline
\textbf{1/3}  & 
$
\begin{aligned}
u_s' &= \hr{0.040 u_s + 0.025 u_i} -0.001 v_i \hr{ -0.021 u_s^2 -0.015 u_s u_i -0.019 u_s v_s -0.029 u_s v_i }\\&+ 0.023 u_i^2 + 0.014 u_i v_s -0.016 u_i v_i -0.001 v_s^2 + 0.002 v_s v_i\\
u_i' &=\hr{  -0.009 u_i + 0.004 u_s u_i + 0.009 u_s v_i -0.014 u_i^2 -0.016 u_i v_s -0.022 u_i v_i}\\
v_s' &= -0.001 u_s + 0.081 u_i \hr{ -0.012 v_s + 0.011 v_i} + 0.006 u_s^2 -0.075 u_s u_i \hr{ + 0.040 u_s v_s} -0.018 u_s v_i \\&-0.118 u_i^2 \hr{+ 0.061 u_i v_s} + 0.001 v_s^2 \hr{-0.011 v_s v_i}\\
v_i' &= -0.004 u_i \hr{-0.020 v_i} -0.001 u_s^2 + 0.020 u_s u_i \hr{+ 0.040 u_s v_i} -0.019 u_i^2 \hr{-0.003 u_i v_s + 0.049 u_i v_i + 0.010 v_s v_i}
\end{aligned}
$
\\  \hline
\textbf{1/2}  & 
$
\begin{aligned}
u_s' &= \hr{0.040 u_s + 0.038 u_i} + 0.001 v_s -0.001 v_i \hr{ -0.020 u_s^2 -0.029 u_s u_i -0.018 u_s v_s -0.035 u_s v_i }\\&+ 0.028 u_i v_s -0.017 u_i v_i -0.001 v_s^2 + 0.002 v_s v_i\\
u_i' &= \hr{ -0.008 u_i + 0.004 u_s u_i + 0.009 u_s v_i -0.015 u_i^2 -0.016 u_i v_s -0.023 u_i v_i}\\
v_s' &= 0.001 u_s -0.031 u_i \hr{ -0.012 v_s + 0.011 v_i} + 0.001 u_s^2 + 0.031 u_s u_i \hr{ + 0.036 u_s v_s} + 0.015 u_s v_i \\&+ 0.066 u_i^2\hr{ -0.034 u_i v_s} + 0.039 u_i v_i + 0.004 v_s^2 \hr{-0.014 v_s v_i} + 0.001 v_i^2\\
v_i' &= 0.013 u_i \hr{-0.020 v_i} + 0.006 u_s u_i + 0.001 u_s v_s \hr{+ 0.037 u_s v_i} -0.045 u_i^2 \hr{+ 0.008 u_i v_s + 0.040 u_i v_i + 0.010 v_s v_i}
\end{aligned}
$
\\  \hline
\textbf{2/3}  & 
$
\begin{aligned}
u_s' &= \hr{0.040 u_s + 0.038 u_i} + 0.001 v_s -0.001 v_i \hr{-0.020 u_s^2 -0.029 u_s u_i -0.019 u_s v_s -0.035 u_s v_i} \\&+ 0.026 u_i v_s -0.017 u_i v_i -0.001 v_s^2 + 0.001 v_s v_i\\
u_i' &= \hr{-0.008 u_i + 0.004 u_s u_i + 0.009 u_s v_i -0.015 u_i^2 -0.016 u_i v_s -0.023 u_i v_i}\\
v_s' &= 0.001 u_s -0.032 u_i \hr{ -0.011 v_s + 0.011 v_i} + 0.001 u_s^2 + 0.032 u_s u_i \hr{ + 0.036 u_s v_s} + 0.015 u_s v_i \\&+ 0.069 u_i^2 \hr{ -0.032 u_i v_s} + 0.039 u_i v_i + 0.003 v_s^2 \hr{ -0.014 v_s v_i} + 0.001 v_i^2\\
v_i' &= 0.012 u_i \hr{-0.020 v_i} + 0.007 u_s u_i + 0.001 u_s v_s \hr{+ 0.037 u_s v_i} -0.044 u_i^2 \hr{ + 0.008 u_i v_s+ 0.041 u_i v_i + 0.010 v_s v_i}
\end{aligned}
$
\\ \hline
\textbf{} & \textbf{Original Equations} \\
\hline
  & 
$
\begin{aligned}
u_s' &= \hr{0.04 u_s + 0.03 u_i} \hr{-0.01 u_s^2 -0.01 u_s u_i -0.02 u_s v_s -0.03 u_s v_i} \\
u_i' &= \hr{-0.01 u_i + 0.01 u_s u_i + 0.01 u_s v_i -0.01 u_i^2 -0.02 u_i v_s -0.02 u_i v_i}\\
v_s' &= \hr{ -0.01 v_s + 0.01 v_i + 0.04 u_s v_s + 0.03 u_i v_s  -0.01 v_s v_i}\\
v_i' &=  \hr{-0.02 v_i + 0.04 u_s v_i + 0.01 u_i v_s+ 0.041 u_i v_i + 0.01 v_s v_i}
\end{aligned}
$
\\  \hline
\end{tabular}
\caption{\rev{SINDy equations identified for the LVSIS temporal system at training proportions $p=1/3$, $1/2$, and $2/3$; the original equations are included for comparison.}}
\label{tab:lvsis}
\end{table}

\begin{figure}[h!]
\centering
\setlength{\tabcolsep}{4pt}
\renewcommand{\arraystretch}{1.2}
\begin{tabular}{c c c c c}
 & {1/3} & {1/2} & {2/3} \\
\centering\rotatebox{90}{NODEs} & 
\includegraphics[width=0.3\textwidth]{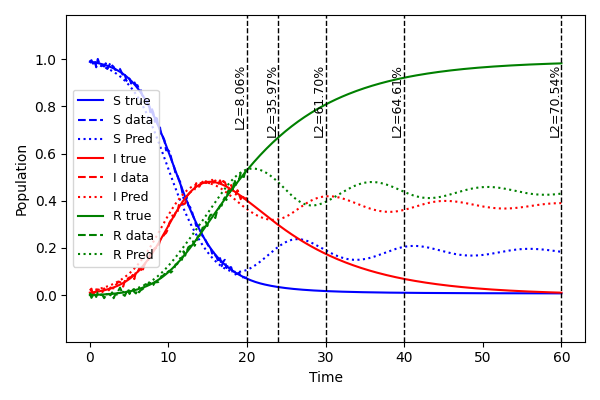} &
\includegraphics[width=0.3\textwidth]{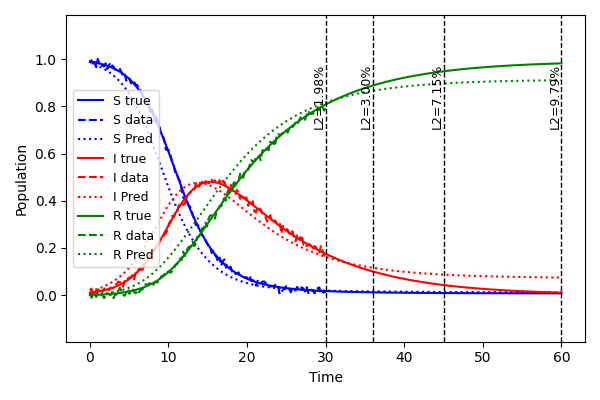} &
\includegraphics[width=0.3\textwidth]{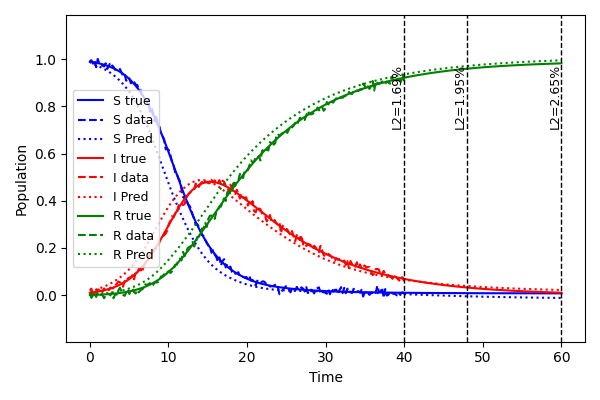} 
\\
\centering\rotatebox{90}{KANODEs} &
\includegraphics[width=0.3\textwidth]{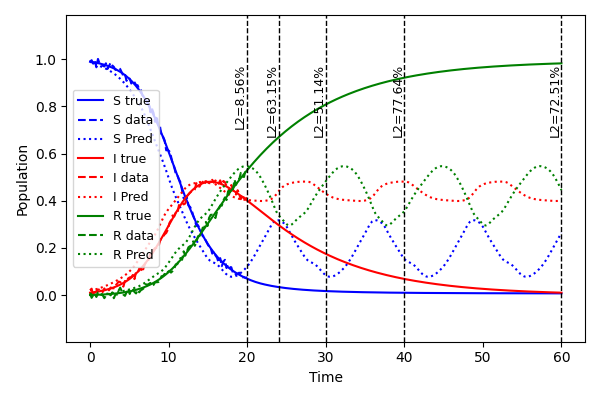} &
\includegraphics[width=0.3\textwidth]{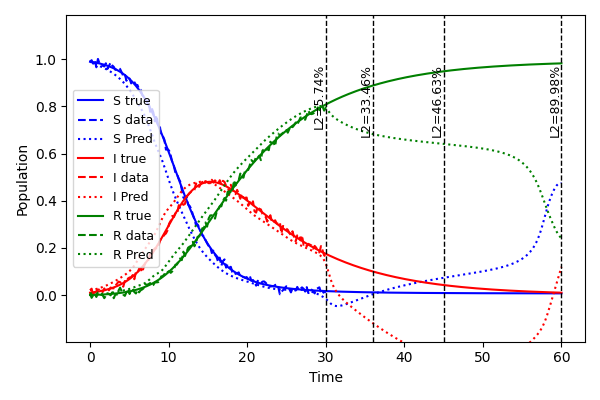} &
\includegraphics[width=0.3\textwidth]{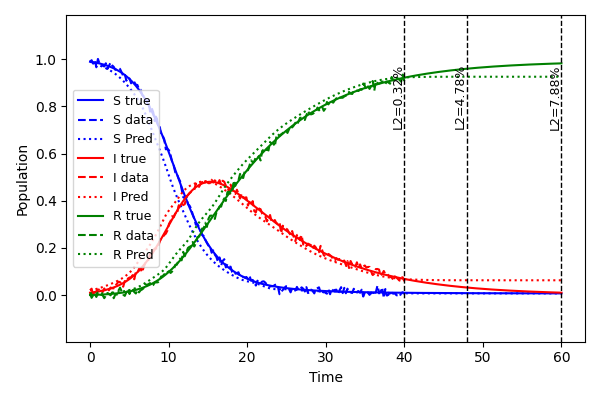} 
\\
\end{tabular}
\caption{\rev{SIR model. Comparison of prediction performance. Noised data.}}
\label{fig:test-sir2}
\end{figure}

\begin{figure}[h!]
\centering
\setlength{\tabcolsep}{4pt}
\renewcommand{\arraystretch}{1.2}
\begin{tabular}{c c c c c}
 & {1/3} & {1/2} & {2/3} \\
\centering\rotatebox{90}{NODEs} & 
\includegraphics[width=0.3\textwidth]{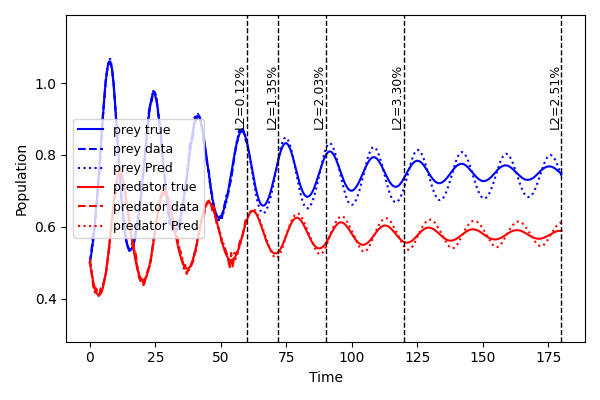} &
\includegraphics[width=0.3\textwidth]{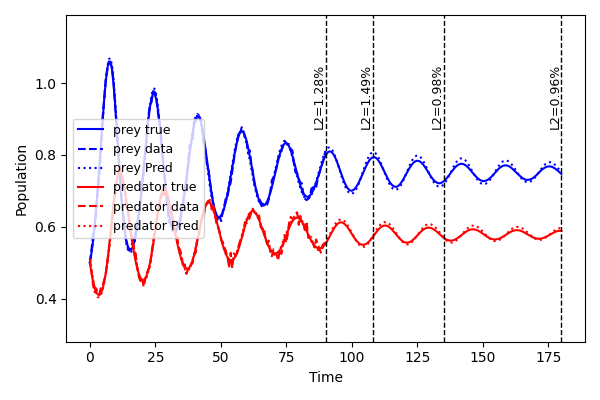} &
\includegraphics[width=0.3\textwidth]{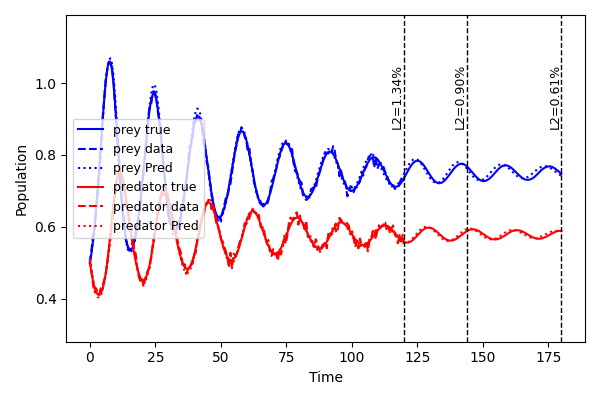} 
\\
\centering\rotatebox{90}{KANODEs} &
\includegraphics[width=0.3\textwidth]{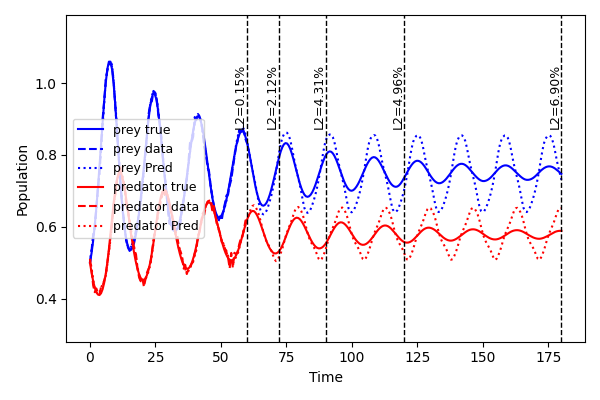} &
\includegraphics[width=0.3\textwidth]{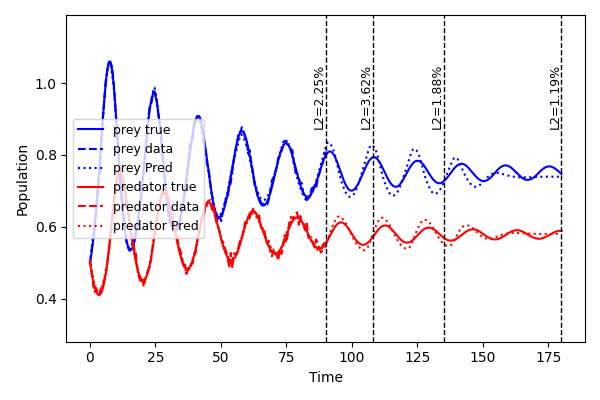} &
\includegraphics[width=0.3\textwidth]{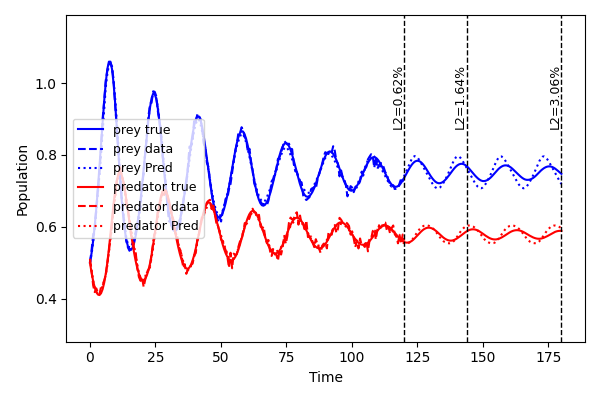} 
\\
\end{tabular}
\caption{\rev{LV model. Comparison of prediction performance. Noised data.}}
\label{fig:test-lv2}
\end{figure}

\begin{figure}[h!]
\centering
\setlength{\tabcolsep}{4pt}
\renewcommand{\arraystretch}{1.2}
\begin{tabular}{c c c c c}
 & {1/3} & {1/2} & {2/3} \\
\centering\rotatebox{90}{NODEs} & 
\includegraphics[width=0.3\textwidth]{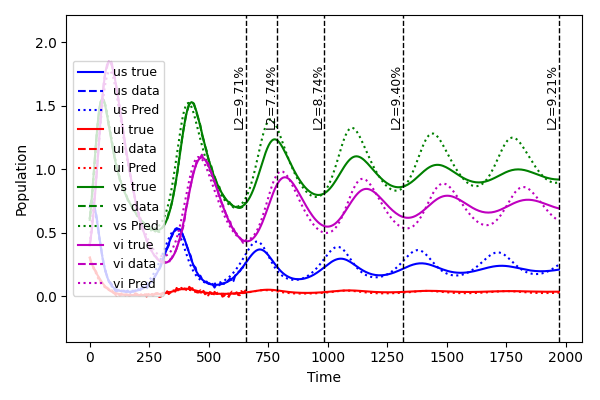} &
\includegraphics[width=0.3\textwidth]{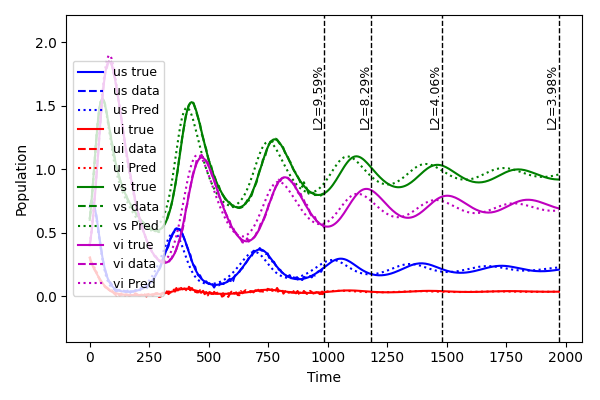} &
\includegraphics[width=0.3\textwidth]{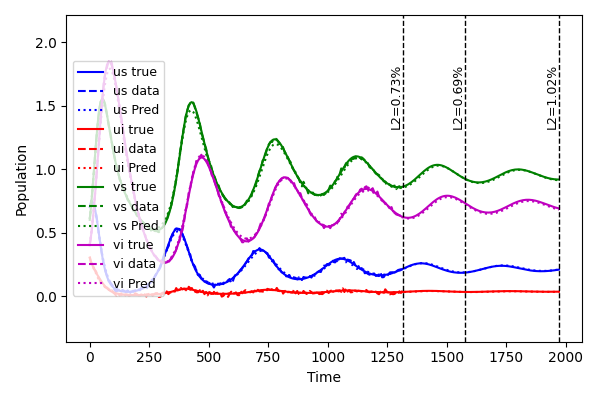} 
\\
\centering\rotatebox{90}{KANODEs} &
\includegraphics[width=0.3\textwidth]{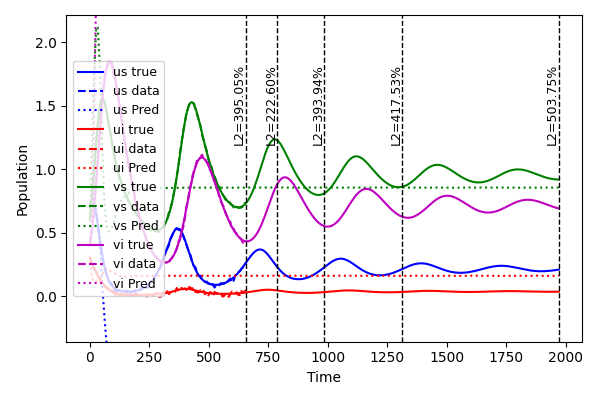} &
\includegraphics[width=0.3\textwidth]{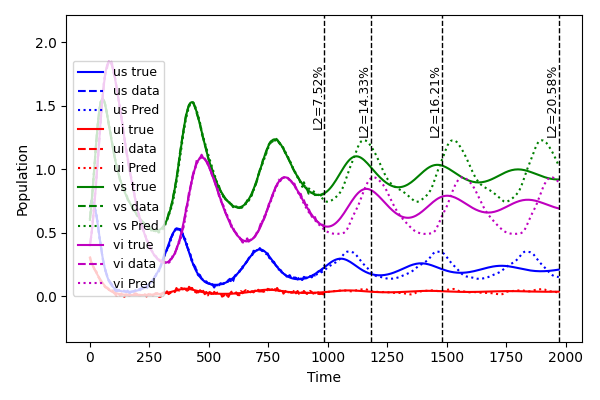} &
\includegraphics[width=0.3\textwidth]{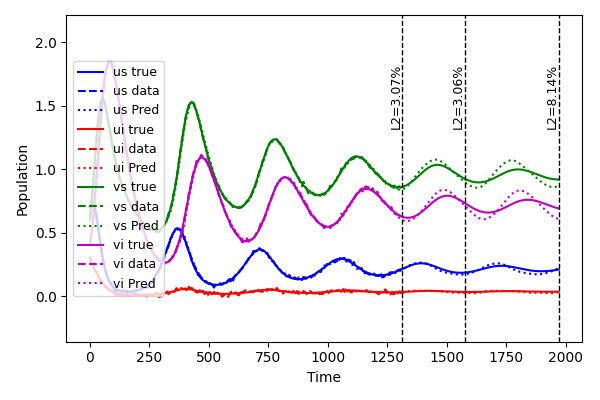} 
\\
\end{tabular}
\caption{\rev{LVSIS model. Comparison of prediction performance. Noised data.}}
\label{fig:test-lvsis2}
\end{figure}

The first row of Figures \ref{fig:test-sir1}, \ref{fig:test-lv1}, and \ref{fig:test-lvsis1} presents the results obtained using the SINDy framework.
The identified ODEs for each model for different training proportions are presented in Tables \ref{tab:sir}, \ref{tab:lv}, and \ref{tab:lvsis}. 
For the SIR model, the learned polynomial terms capture the infection and recovery dynamics. As the training set increases from $1/3$ to $2/3$, the coefficients approach the true parameters, improving accuracy and stability.  
In the LV system, the learned equations reproduce the characteristic predator–prey oscillations. The bilinear interaction $uv$ is consistently identified, while quadratic terms $u^2$ and $v^2$ provide small corrections for nonlinear effects. Larger training sets improve the match to reference trajectories, particularly in amplitude and period.  
For the LVSIS eco-epidemiological model, SINDy recovers both intra- and inter-species interactions, including bilinear transmission and predation terms and nonlinear density-dependent effects. Even small training sets capture the qualitative dynamics, while larger sets yield more accurate coefficients and closer agreement with reference trajectories. 
Overall, we observe that the SINDy framework is capable of capturing the essential structure of the underlying dynamics across all three systems.

In the second and third rows of Figures \ref{fig:test-sir1}, \ref{fig:test-lv1}, and \ref{fig:test-lvsis1}, we present the results obtained using the two neural network approaches, NODEs and KANODEs. Both methods achieve a good performance, providing accurate reconstructions of the underlying trajectories and demonstrating robust predictive capabilities, particularly with larger training datasets. Increasing the amount of training data consistently enhances prediction accuracy for all methods over all datasets. However, for the SIR dataset, the model exhibits reduced accuracy when predicting beyond the training set.  For the more complex oscillatory dynamics of the LV and LVSIS models, the models achieve good predictive performance. 
KANODEs perform well on larger datasets ($p=2/3$). Neural ODEs yield better extrapolation outside the training set for all models. In the Lotka--Volterra system, both methods achieve good performance across all training proportions, reflecting the richer dynamics. For the LVSIS model, good results are observed for NODEs with $p=1/2$ and for KANODEs with $p=2/3$. We note that, for Neural ODEs, we present results using a larger network (NODEs$_2$), which contains significantly more parameters than KANODEs.

\rev{Next, we present the results obtained for noisy data with a noise level of $\sigma = 0.01$. Figures~\ref{fig:test-sir2}, \ref{fig:test-lv2}, and \ref{fig:test-lvsis2} show the corresponding results obtained using the neural network approaches.}
The SIR model, which exhibits relatively simple dynamics, is highly sensitive to the size of the training set. However, Neural ODEs provide accurate predictions when trained with sufficient temporal information, particularly for $p = 2/3$.  
The richer and more oscillatory trajectories in the LV dataset allow the models to effectively learn the underlying dynamics, even from smaller training sets. We observe that Neural ODEs generally outperform KANODEs, especially when trained on limited data, by better capturing nonlinear interactions and coupled oscillatory behavior among the species. Overall, both approaches achieve low prediction error when provided with representative training data that adequately characterizes the dynamics of the underlying system.

\section{\rev{Temporal Wall-Clock Results}}
\label{app:temporal-timing}

\rev{In this appendix the timing evidence separately from the main accuracy tables. All wall-clock measurements were obtained on the hardware described in Section~4. Neural  network models were trained for $500$ epochs using the Adam optimizer with learning rate $5\times10^{-4}$. Each epoch used $20$ optimization steps for SIR and LV and $10$ optimization steps for LVSIS.  }

\begin{table}[h!]
\centering
\small
\begin{tabular}{|l|c|c|c|}
\hline
 & 1/3 & 1/2 & 2/3 \\
\hline
\multicolumn{4}{|c|}{SIR} \\
\hline
SINDy & 0.010 & 0.003 & 0.002 \\
NODE$_1$ & 32.06 $\pm$ 0.08 & 32.44 $\pm$ 0.45 & 32.53 $\pm$ 0.06 \\
NODE$_2$ & 82.22 $\pm$ 0.27 & 80.08 $\pm$ 2.45 & 83.12 $\pm$ 0.47 \\
KANODE & 122.5 $\pm$ 0.88 & 127.2 $\pm$ 4.23 & 121.6 $\pm$ 1.03 \\
\hline
\multicolumn{4}{|c|}{LV} \\
\hline
SINDy & 0.002 & 0.002 & 0.002 \\
NODE$_1$ & 30.82 $\pm$ 0.59 & 30.49 $\pm$ 0.24 & 31.08 $\pm$ 0.44 \\
NODE$_2$ & 79.86 $\pm$ 0.92 & 79.84 $\pm$ 0.40 & 79.73 $\pm$ 0.33 \\
KANODE & 110.2 $\pm$ 0.73 & 110.7 $\pm$ 0.30 & 110.7 $\pm$ 0.54 \\
\hline
\multicolumn{4}{|c|}{LVSIS} \\
\hline
SINDy & 0.006 & 0.003 & 0.003 \\
NODE$_1$ & 27.74 $\pm$ 0.37 & 27.54 $\pm$ 0.10 & 28.15 $\pm$ 0.04 \\
NODE$_2$ & 77.49 $\pm$ 1.08 & 76.12 $\pm$ 0.31 & 78.77 $\pm$ 0.40 \\
KANODE & 128.5 $\pm$ 1.19 & 129.2 $\pm$ 0.92 & 131.2 $\pm$ 1.17 \\
\hline
\end{tabular}
\caption{\rev{Wall-clock fitting time (s) for SINDy and training time (s) for the neural models on noise-free temporal data at training proportions $p=1/3$, $1/2$, and $2/3$.}}
\label{tab:temporal-training-time}
\end{table}

\begin{table}[h!]
\centering
\small
\begin{tabular}{|l|c|c|c|c|}
\hline
\rev{Epoch} & Train & Test & Full & Time (s) \\
\hline
\multicolumn{5}{|c|}{NODE$_1$} \\
\hline 
 25 & 17.1575 & 7.8158 & 12.4866 & 1.37  \\
 50 & 12.2175 & 8.0906 & 10.1540 & 2.73  \\
 100 & 8.3540 & 11.9554 & 10.1547 & 5.43 \\
 250 & 5.5846 & 4.8469 & 5.2158 & 13.55  \\
 500 & 4.1408 & 2.8319 & 3.4863 & 27.08  \\
 1000 & 2.4789 & 2.1065 & 2.2927 & 54.18  \\
 2000 & 0.7786 & 0.8211 & 0.7999 & 108.47  \\
 4000 & 0.6627 & 0.6199 & 0.6413 & 217.07  \\
 6000 & 0.7760 & 0.7106 & 0.7433 & 328.09  \\
\hline 
\multicolumn{5}{|c|}{NODE$_2$} \\
\hline 
 25 & 14.1723 & 8.1877 & 11.1800 & 3.83  \\
 50 & 9.0532 & 12.6916 & 10.8724 & 7.64  \\
 100 & 12.7763 & 34.4451 & 23.6107 & 15.62  \\
 250 & 2.6151 & 4.8373 & 3.7262 & 38.51  \\
 500 & 2.6762 & 1.3634 & 2.0198 & 77.50  \\
 1000 & 0.7940 & 0.9539 & 0.8740 & 155.90  \\
 2000 & 0.5170 & 0.8113 & 0.6641 & 311.75  \\
 4000 & 0.7754 & 0.7062 & 0.7408 & 619.54  \\
 6000 & 0.9612 & 0.6697 & 0.8154 & 929.66  \\
\hline
\multicolumn{5}{|c|}{KANODE} \\
\hline 
 25 & 15.8686 & 14.7313 & 15.2999 & 6.88  \\
 50 & 4.9209 & 9.3456 & 7.1333 & 13.77  \\
 100 & 1.2087 & 6.7127 & 3.9607 & 27.23  \\
 250 & 0.6081 & 3.0820 & 1.8450 & 66.61  \\
 500 & 0.2307 & 2.2879 & 1.2593 & 132.12  \\
 1000 & 0.1714 & 1.9722 & 1.0718 & 260.35  \\
 2000 & 0.1049 & 1.5483 & 0.8266 & 518.79  \\
\hline
\end{tabular}
\caption{\rev{Training-, test-, and full-trajectory relative mean errors (\%) and cumulative training time (s) for LVSIS at training proportion $p=1/2$.}}
\label{tab:temporal-convergence}
\end{table}

\rev{Table~\ref{tab:temporal-training-time} provides a wall-clock comparison for the original temporal data. The neural models were trained for $500$ epochs using Adam with learning rate $5\times10^{-4}$; each epoch used $20$ optimization steps for SIR and LV and $10$ for LVSIS. Neural-model values are reported as mean $\pm$ standard deviation over random initializations, whereas deterministic SINDy fits are reported as single values. Direct SINDy fitting takes $0.002$--$0.010$~s because it solves a sparse regression problem after the candidate library and numerical derivatives have been constructed. The neural approaches require repeated differentiable ODE rollouts: NODE$_1$ takes approximately $28$--$33$~s, NODE$_2$ takes $76$--$83$~s, and KANODE takes $110$--$131$~s. Although KANODE has substantially fewer trainable parameters, it is about $3.6$--$4.7$ times slower than NODE$_1$ and $1.4$--$1.7$ times slower than NODE$_2$ in these experiments. The spline evaluations and their derivatives therefore make the parameter count an incomplete measure of computational cost. SINDy's substantially shorter fitting time should be interpreted with caution, as it relies on a prescribed candidate library and accurate derivative estimates. Moreover, the same direct fitting procedure did not produce stable models for noisy data.}

\rev{Table~\ref{tab:temporal-convergence} complements the aggregate timing comparison by showing how the LVSIS errors evolve with both the number of training epochs and cumulative wall-clock time. Reported errors correspond to the model with the lowest training loss achieved up to each epoch. 
KANODE is more expensive per epoch: at 2000 epochs it requires $518.79$~s, compared with $108.47$~s for NODE$_1$ and $311.75$~s for NODE$_2$. However, its error decreases more rapidly during the early part of training, so fewer epochs can be sufficient for a moderate target accuracy. For example, after only $50$ epochs, KANODE has a full-horizon error of $7.13\%$, compared with $10.15\%$ for NODE$_1$ and $10.87\%$ for NODE$_2$. KANODE reduces the error below $2\%$ after $250$ epochs, reaching $1.85\%$ in $66.61$~s; NODE$_2$ first falls below $2\%$ after $1000$ epochs, reaching $0.87\%$ in $155.90$~s, while NODE$_1$ first crosses this threshold after $2000$ epochs, reaching $0.80\%$ in $108.47$~s. Thus, for this accuracy threshold, KANODE uses four to eight times fewer epochs and also less total wall-clock time despite its higher per-epoch cost. These results show that computational efficiency depends on both the cost per epoch and the number of epochs required to attain the desired error.}

\end{document}